\documentstyle[11pt]{article}
\begin{document}

\newcommand{\ul}{\underline}
\newcommand{\sa}{\Sigma}

\title{ TOPOLOGICAL QUANTUM FIELD THEORY WITH CORNERS  BASED ON THE KAUFFMAN
BRACKET}
\author{ R{\u{a}}zvan Gelca}
\maketitle

\begin{abstract}
We describe the construction of a topological quantum field theory with corners
based on the Kauffman bracket, that underlies the smooth theory
of Lickorish. We also exhibit some properties of invariants of
3-manifolds with boundary.
\end{abstract}

\begin{center}
{\bf 1. INTRODUCTION}
\end{center}

\bigskip

The discovery of the Jones polynomial [J] and its placement in the 
context of quantum field theory [Wi] led to various constructions of
topological quantum field theories (shortly TQFT's). A first example
of such a theory has been produced by Reshetikhin and Turaev in [RT], and
makes use of the representation theory of Hopf algebras. Then, an
alternative construction based on geometric techniques has been worked out
by Kohno in [Ko]. A combinatorial approach based on skein spaces associated
to the Kauffman bracket [K1] has been exhibited by Lickorish in [L1]
and [L2] and by Blanchet, Habegger, Masbaum, and Vogel in [BHMV].
All these theories are smooth, in the sense that gluings occur along 
closed surfaces, so they are based on the Atiyah-Segal set of axioms [A].

In an attempt to give a more axiomatic approach to such a theory, and also to
make it easier to handle, K. Walker described in [Wa] a system of axioms
for a TQFT in which one allows gluings along surfaces with boundary, a
so called TQFT with corners. He also described the minimal amount of initial
information (basic data) that one needs to know in order to be able to
recover the whole theory from axioms. He based his theory on the decompositions
of surfaces into disks, annuli and pairs of pants, and along with the mapping
class group of a surface he considered the groupoid of transformations 
of these decompositions.

Following partial work from [Wa], in [FK] and [G1] a TQFT
with corners associated to the Reshetikhin-Turaev 
theory has been exhibited . We mention that
in this construction one encounters a sign problem at the level of the 
groupoid of transformations of decompositions.  The presence of this sign
problem was due to the fact [G2] that the theory was based on the Jones
polynomial, whose skein relation is defined for oriented links.

In the present paper we intend to construct a TQFT with corners that underlies
the smooth TQFT of Lickorish [L1], [L2]. It is based on the skein theory of the
Kauffman bracket. Note that since the Kauffman bracket is defined for 
unoriented links, we will not encounter any sign problem this time.
The main elements involved in our construction are the Jones-Wenzl
idempotents [We]. They are the analogues of the irreducible representations
from [RT]. Recall the cabling formula from [KM], which enables us to
compute (global) link invariants via skein formulas. The Jones-Wenzl 
idempotents help us do  this kind of computations locally. Regarding 
the computations, we make the observation that in our case they will be done
either in the skein space of the plane, or in that of the disk with points on
the boundary, although 
the spaces associated to closed surfaces are skein spaces 
of handlebodies [L2], [R].

The paper has many similarities with [FK] and [G1], but from many points
of view it is much simpler. For example, we no longer need a special
direction in the plane, like in the case of computations with representations, 
thus the problem of rotating trivalent vertices is trivial.
Also, the analogues of the Clebsch-Gordan coefficients are easier to
construct and comprehend.

In Section 2 we review the definitions from [Wa]. Section 3 starts with
a review of facts about skein spaces and then proceeds with the description
of the basic data. In Section 4 we prove 
that the basic data gives rise to a well
defined TQFT. As a main device involved in the proof we exhibit
a tensor contraction formula. In the last section we  generalize 
to surfaces with boundary a well known formula for the invariant of the
product of a closed surface with a circle. Next we show that the invariants
of 3-manifolds with boundary have a distinguished vector component
which satisfies the Kauffman bracket skein relation. As a consequence,
we compute the invariant of the complement of a regular neighborhood
of a link, and explain how the invariants of closed manifolds arise when
doing surgery on such links.

\begin{center}
{\bf 2. FACTS ABOUT TQFT'S WITH CORNERS}
\end{center}

A TQFT with corners is one that allows gluings of 3-manifolds
along surfaces in 
their boundary that themselves have boundary. In order to be able to
understand such a theory we must first briefly describe its objects, 
the extended surfaces and 3-manifolds. For an extensive discussion
we recommend [Wa]. The adjective ``extended'' comes from the way 
 the projective ambiguity of the invariants is resolved, which is done, as
usually, via an extension of the mapping class group. All surfaces and
3-manifolds throughout the paper are supposed to be piecewise linear, compact
and orientable.

In order to fulfill the needs of a TQFT with corners, the concept of extended
surface will involve slightly more structure than the usual Lagrangian space,
namely the decomposition into disks, annuli and pairs of pants (shortly
DAP-decomposition).

\medskip

\ul{Definition.} A DAP-decomposition of a surface $\Sigma $ consists of

\ \ - a collection of disjoint simple closed curves in the interior of
$\Sigma $ that cut $\Sigma $ into elementary surfaces: disks, annuli, and 
pairs of pants, and an ordering of these elementary surfaces;

\ \ - a numbering of the boundary components of each elementary surface 
 $\Sigma _0$ by 1 if $\sa _0$ is a disk, 1 and 2 
if $\sa _0$ is an annulus, and 1, 2 and 3 if $\sa _0$ is a pair of pants;

\ \ - a parametrization of each boundary component $C$ of 
$\sa _0$ 
by $S^1=\{z|\ |z|=1\}$ (the parametrization being compatible with the
orientation of $\sa _0$ under the convention ``first out'');

\ \ -  fixed disjoint embedded arcs
in $\sa _0$ joining $e^{i\epsilon}$ (where  $\epsilon >0$ is small)
on the $j$-th boundary component to
$e^{-i\epsilon}$ on the $j+1$-st (modulo the number of boundary
components of $\sa _0$). These arcs will be called seams.

\medskip

An example of a DAP-decomposition is shown in Fig 2.1. Two decompositions are
considered identical if they coincide up to isotopy. We also make the
convention that whenever we talk about the decomposition curves we also
include the boundary components of the surface as well.

\input epsf

\begin{figure}[htbp]
\centering
\leavevmode
\epsfxsize=3.5in
\epsfysize=1.1in
\epsfbox{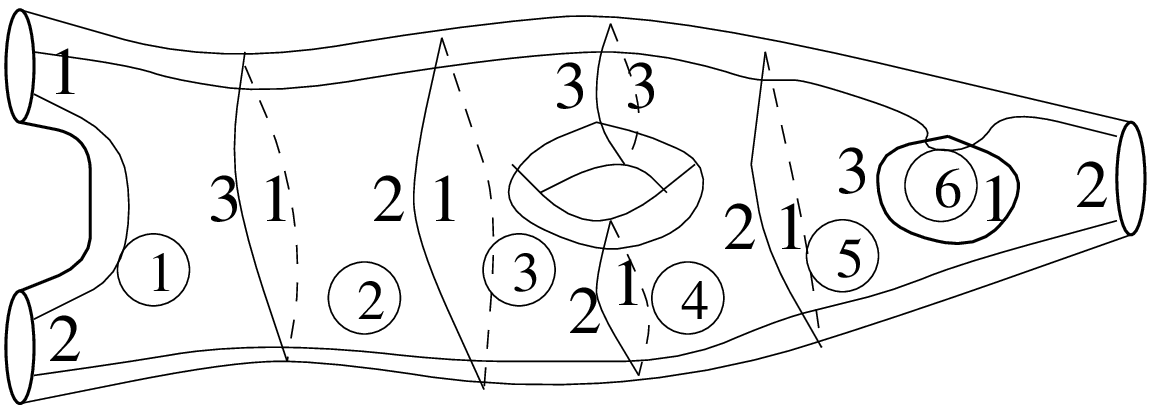}

Fig. 2.1.
\end{figure}

\medskip

\ul{Definition.} An extended surface (abbreviated e-surface) is a pair
$(\sa ,D)$ where $\sa $ is a surface and $D$ is a DAP-decomposition 
of $\sa $.

\medskip

Let us note that in the case of smooth TQFT's one is only interested
in the Lagrangian subspace spanned by the decomposition curves of $D$
in $H_1(\sa )$. 
In our case, we will be interested in the decomposition itself, since 
we can always arrange the gluing to be along a collection of elementary
subsurfaces in the boundary of the 3-manifold. We emphasize that the
DAP-decomposition plays the same role as the basis plays for a vector 
space. 

If we change the orientation of a surface, the DAP-decomposition should
be changed by reversing all orientations and subsequently by permuting the 
numbers 2 and 3 in the pairs of pants.

In what follows, we will call a move any transformation of one 
DAP-decomposition into another. By using Cerf theory [C] one can 
show that any move can be written as a composition of the elementary moves
described in Fig. 2.2 and their inverses, together with the permutation map
$P$ that changes the order of elementary surfaces. In the 
sequel $T_1$ will be called a twist, $R$ rotation, the maps $A$ and
$D$ contractions of annuli, respectively disks, and their inverses
expansions of annuli and disks.

\begin{figure}[htbp]
\centering
\leavevmode
\epsfxsize=5.2in
\epsfysize=3.9in
\epsfbox{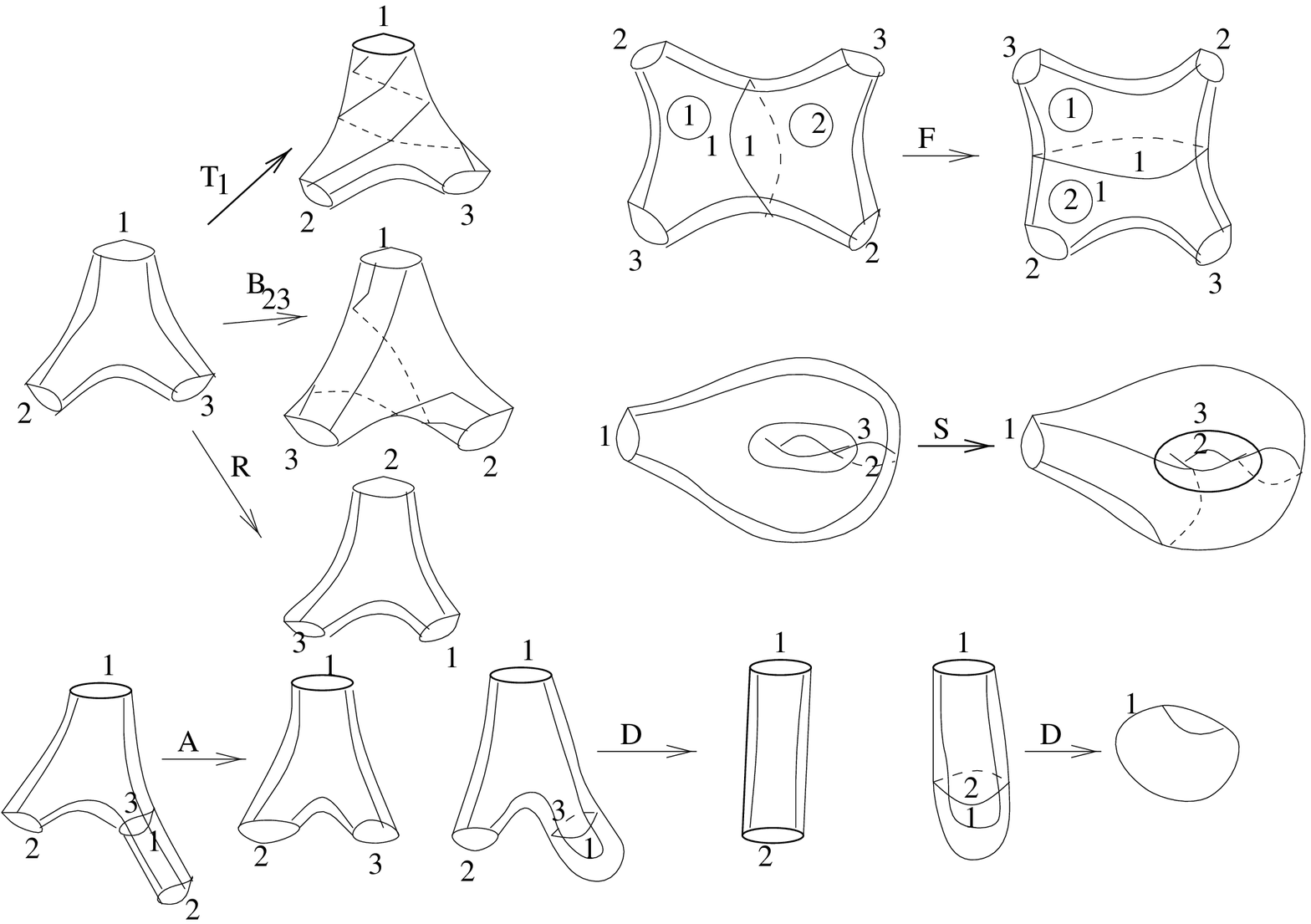}

Fig. 2.2.
\end{figure}

\medskip

\ul{Definition.} An extended morphism (shortly e-morphism) is a map
between two e-surfaces
$(f,n):(\sa _1,D_1)\rightarrow (\sa _2,D_2)$ where $f$ is a homeomorphism
and $n$ is an integer.

\medskip

Note that such an e-morphism can be written as a composition of
a homeomorphism $(f,0):((\sa _1,D_1)\rightarrow (\sa _2,f(D_1))$, a move
$(\sa _2,f(D_1))\rightarrow (\sa _2,D_2)$
and the morphism $(0,n):(\sa _2, D_2)\rightarrow (\sa _2, D_2)$.
Note also that the moves from Fig 2.2 have the associated homeomorphism
equal to the identity.

The set of e-morphisms is given a groupoid structure by means of the following
composition law.
For $(f_1,n_1):(\sa _1,D_1)\rightarrow (\sa _2,D_2)$ and $(f_2,n_2):
(\sa _2, D_2)\rightarrow ( \sa _3,D_3)$ let

\ \ $(f_2,n_2)(f_1,n_1):=(f_2f_1,n_1+n_2-\sigma ((f_2f_1)_*L_1,(f_2)_*L_2
,L_3)$

\noindent where $\sigma $ is Wall's nonadditivity function [W] and $L_i
\subset H_1(\sa _i)$ is the subspace generated by the decomposition 
curves of $D_i$, $i=1,2,3$.

Let us now review some facts about extended 3-manifolds.

\medskip

\ul{Definition.} The triple $(M,D,n)$ is called an extended 3-manifold 
(e-3-manifold) if $M$ is a 3-manifold, $D$ is a DAP-decomposition
of $\partial M$ and $n\in {\bf Z}$.

\medskip

The boundary operator, disjoint union and mapping cylinder are defined in
the canonical way, namely $\partial (M,D,n)=(\partial M, D)$,
$(M_1,D_1,n_1)\sqcup (M_2,D_2,n_2)=(M_1\sqcup M_2, D_1\sqcup D_2, n_1+n_2)$
and for $(f,n):(\sa _1,D_1)\rightarrow (\sa _2,D_2)$, $I_{(f,n)}=(I_f,D,n)$,
with the only modification that in $I_f$ we identify the boundary components
of $- \sa _1$ with those of $\sa _2$ that they get mapped onto, thus $\partial
I_f= -\sa _1 \cup \sa _2$ and $D=D_1\cup D_2$. More complicated is the
gluing of e-3-manifolds, which is done as follows.

\medskip

\ul{Definition.} Let $(M,D,n)$ be an e-3-manifold and $(\sa _1, D_1)$
and $(\sa _2,D_2)$ be two disjoint surfaces in its boundary. Let
$(f,m):(\sa _1,D_1)\rightarrow (\sa _2,D_2)$ be an e-morphism.
Define the gluing of $(M,D,n)$ by $(f,m)$ to be

\ \ \ $(M,D,n)_{(f,m)}:=(M_f,D', m+n-\sigma (K,L_1\oplus L_2,\Delta^-))$

\noindent where $M_f$ is the gluing of $M$ by $f$, $D'$ is the 
image of $D$ under this gluing, $\sigma $ is Wall's nonadditivity
function, $K$ is the subspace of $H_1(\partial M)$ spanned
by the kernel of $H_1(\sa _1\cup \sa _2)\rightarrow H_1(M)$, $\partial
\sa _1$ and $\partial \sa _2$, $L_i$ are the subspaces of $H_1(\sa _i)$
generated by the decomposition curves of $D_i$ and $\Delta ^-=\{(x,-f_*(x)),
x\in H_1(\sa _1)\}$.

\medskip

For a geometric explanation of this definition see [Wa].

In order to define a TQFT we also need a finite set of labels ${\cal L}$,
with a distinguished element $0\in {\cal L}$. Consider the category
of labeled extended surfaces (le-surfaces) whose objects 
are e-surfaces with the 
boundary components labeled by elements in ${\cal L}$ (le-surfaces),
 and whose morphisms
are the e-morphisms that preserve labeling (called labeled extended 
morphism and abbreviated le-morphisms).
An le-surface is thus a triple $(\sa ,D,l)$, where $l$ is a labeling
function.

Following [Wa] we define a TQFT with label set ${\cal L}$ to consist out
of 

-a functor $V$ from the category of le-surfaces to that of finite dimensional
vector spaces, called modular functor,

-a partition function $Z$ that associates to each 3-manifold a vector in the
vector space of its boundary.

The two should satisfy the following axioms:

(2.1) (disjoint union) $V(\sa _1\sqcup \sa _2,D_1\sqcup D_2,l_1\sqcup l_2)=
V(\sa _1,D_1,l_1)\bigotimes V(\sa _2,D_2,l_2)$;

(2.2) (gluing for $V$) Let $(\sa ,D)$ be an le-surface, $C,C'$ two subsets 
of boundary components of $(\sa )$, and $g:C\rightarrow C'$ the homeomorphism
which is the parametrization reflecting map. Let $\sa _g$ be the gluing
of $\sa $ by $g$, and $D_g$ the DAP-decomposition induced by $D$. Then, for
a certain labeling $l$ of $\partial \sa $ we have

\ \ \ \ \ \ \ 
\ \ \ $V(\sa _g, D_g, l)=\bigoplus _{x\in {\cal L}(C)}V(\sa ,D,(l,x,x))$

\noindent where the sum is over all labelings of $C$ and $C'$ by $x$.

(2.3) (duality) $V(\sa ,D,l)^*=V(-\sa, -D,l)$  
and the identifications
$V(\sa,D ,l)=V(-\sa, -D, l)^*$ and $V(-\sa, -D, l)=V(\sa, D ,l)^*$
are mutually adjoint. Moreover, the following conditions should be satisfied

-if $(f,n)$ is an le-morphism between to le-surfaces,
 then $V(\bar{f},-n)$ is the adjoint inverse of $V(f,n)$, where we
denote by $\bar{f}$ the homeomorphism
induced between the surfaces with reversed orientation.

-if $\alpha _1\otimes \alpha _2\in V(\sa _1,D_1,l_1)\bigotimes 
V(\sa _2,D_2,l_2)$ and $\beta _1\otimes \beta _2\in
V(-\sa _1,-D_1,l_1)\bigotimes 
V(-\sa _2,-D_2,l_2)$
then $<\alpha _1\otimes \alpha _2,\beta _1\otimes \beta _2>=<\alpha _1,\beta 
_1>
<\alpha _2,\beta _2>$,

-there exists a function $S:{\cal L}\rightarrow {\bf C}^*$ such that
with the notations from axiom (2.2) if $\oplus _x \alpha _x\in \bigoplus _{
x\in {\cal L}(C)} V(\sa ,D,(l,x,x))$ and $\oplus _x \beta  _x\in \bigoplus _{
x\in {\cal L}(C)} V(-\sa ,-D,(l,x,x))$
then the pairing on the glued surface is given by $<\oplus _x \alpha _x,
\oplus _x \beta  _x>=\sum_x S(x)<\alpha _x,\beta _x>$, where $x=(x_1,x_2,
\cdots ,x_n)$ and $S(x)=S(x_1)S(x_2)\cdots S(x_n)$;

(2.4) (empty surface) $V(\emptyset )={\bf C}$;

(2.5) (disk) If ${\bf D}$ is a disk $V({\bf D},m)={\bf C}$ 
if $m=0$ and $0$ otherwise;

(2.6) (annulus) If $A$ is an annulus then $V(A,(m,n))={\bf C}$ if
$m=n$ and $0$ otherwise;

(2.7) (disjoint union for $Z$) $Z((M_1,D_1,n_1)\sqcup (M_2,D_2,n_2))=
Z(M_1,D_1,n_1)\otimes Z(M_2,D_2,n_2)$;

(2.8) (naturality) Let $(f,0):(M_1,D,n)\rightarrow (M_2,f(D),n)$. Then
$V(f|\partial (M_1,D,n))Z(M_1,D,n)=Z(M_2,f(D),n)$.

(2.9) (gluing for $Z$) Let $(\sa _1,D_1)$, $(\sa _2, D_2)\subset 
\partial (M,D,m)$ be disjoint, and let 
$(f,n):(\sa _1,D_1)\rightarrow(\sa _2, D_2) $.
Then by (2.2) 

\noindent $V(\partial (M,D,m))=\bigoplus _{l_1,l_2}
V(\sa _1,D_1,l_1)\otimes V(\sa _2,D_2,l_2)\otimes
V(\partial (M,D,m)\backslash ((\sa _1,D_1)\cup (\sa _2,D_2), l_1\cup l_2)$
 
\noindent hence $Z(M,D,m)= \bigoplus_{l_1,l_2}\sum_j \alpha _{l_1}^{(j)}
\otimes \beta _{l_2}^{(j)}\otimes \gamma _{l_1,l_2}^{(j)}$. 
The axiom states that

\ \ \ \ \ \ \ 
$Z((M,D,m)_{(f,n)})=\oplus _l \sum_j<V(f,n)\alpha 
^{(j)}_l,\beta ^{(j)}_l>\gamma
^{(j)}_{ll}$,

-\noindent where $l$ runs through all labelings of $\partial \sa _1$;

(2.10) (mapping cylinder axiom) For $(id,0):(\sa ,D)\rightarrow (\sa ,D)$
we have 

\ \ \ \ \ \ \ \ \ \ \  $Z(I_{(id,0)})=\oplus 
_{l\in {\cal L}(\partial \sa)}id_l$

\noindent where $id_l$ is the identity matrix in $V(\sa ,D,l)\bigotimes
V(\sa ,D,l)^*=V(\sa ,D,l)\bigotimes
V(-\sa ,-D,l)$.

\bigskip

\begin{center}
{\bf 3. THE BASIC DATA}
\end{center}

\bigskip

In order to construct a TQFT with corners one needs to specify a certain
amount of information, called basic data, from which the modular functor and
partition function can be recovered via the axioms. Note that the 
partition function is completely determined by the modular
functor, so we only need to know that latter. Moreover, the modular
functor is determined by the vector spaces associated to le-disks, annuli
and pairs of pants, and by the linear maps associated to le-morphisms.
An important observation is that the matrix of a morphism $V(f,0)$,
where $(f,0):(\sa _1,D)\rightarrow (\sa _2,f(D))$, is the identity matrix, so 
one only needs to know the values of the functor for moves, hence for the
elementary moves described in Fig. 2.2. Of course we also need to know
its value for the map $C=(id,1)$.

The possibility of relating our theory to the Kauffman bracket
depends on the choice of basic data.  Our construction is
inspired by [L1]. We will review the notions we need from that paper
and then proceed with our definitions.

Let $\sa $ be a surface with a collection of $2n$ points on 
its boundary $(n\geq
0)$. A link diagram in $\sa $ is an immersed compact 1-manifold $L$ 
in $\sa $
with the property that $L\cap \partial \sa =\partial L$, $\partial L $ consists
of the $2n$ distinguished point on $\partial \sa $, the singular points
of $L$ are in the interior of $\sa $ and are transverse double points,
and for each such point the ``under'' and ``over'' information is recorded.

Let $A\in {\bf C}$ be fixed. The skein vector space of $\sa $, denoted
by ${\cal S}(\sa )$, is defined to be the complex vector space spanned by
all link diagrams factored by the following two relations:

\ \ a). $L\cup $(trivial closed curve)$=-(A^2+A^{-2})L$,

\ \ b). $L_1=AL_2+A^{-1}L_3$,

\noindent where $L_1,L_2$ and $L_3$ are any three diagrams that coincide 
except in a small disk, where they look like in Fig. 3.1.

\begin{figure}[htbp]
\centering
\leavevmode
\epsfxsize=2.2in
\epsfysize=0.7in
\epsfbox{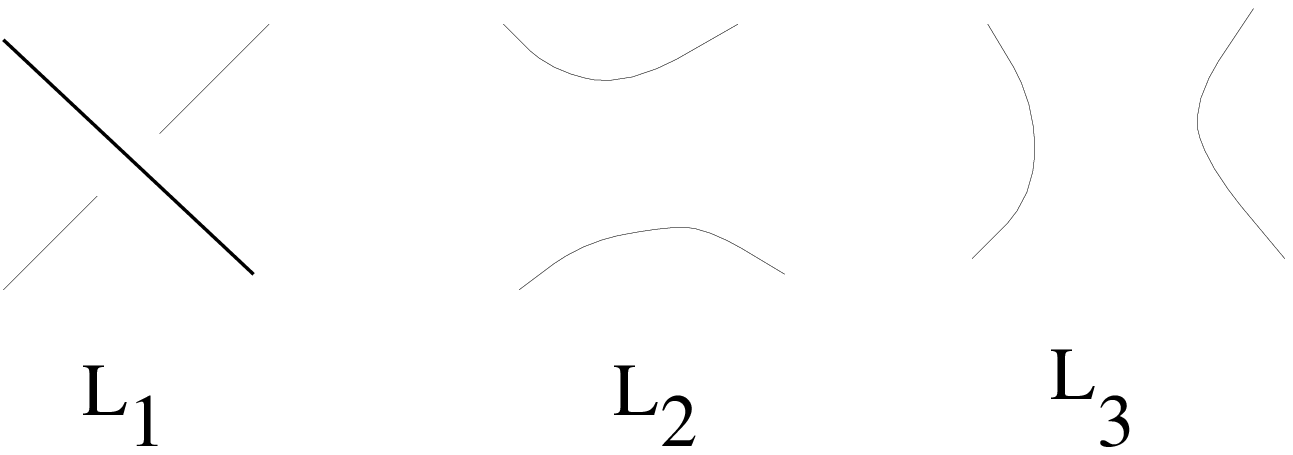}

Fig. 3.1.
\end{figure}

For simplicity, from now on, whenever in a diagram we have an integer, say $k$,
written next to a strand we will actually mean that we have $k$ parallel
strands there.  Also rectangles (coupons) inserted in diagrams will stand
for elements of the skein space of the rectangle inserted there.

Three examples are useful to consider. The first one is the skein space of the
plane, which is the same as the one of the sphere, and it is well known that
it is isomorphic to ${\bf C}$.  

The second example is that of an annulus $A$
with no points on the boundary. It is also a well known fact that
${\cal S}(A)$ is isomorphic to the ring of polynomials ${\bf C}[\alpha ]$,
(if  endowed  with the multiplication defined by the gluing of annuli).
The independent variable $\alpha$ is the diagram with one strand parallel
to the boundary of the annulus. 
Recall from [L1] that every link diagram $L$ in the plane determines a map 

\ \ \ \ \ \ 
\ \ \ \ $<\cdot,\cdot ,\cdots ,\cdot>_L:
{\cal S}(A)\times :{\cal S}(A)\times\cdots :{\cal S}(A)
\rightarrow {\cal S}({\bf R}^2)$

\noindent obtained by first expanding each component of $L$ to an
annulus via the blackboard framing and then homeomorphically mapping $A$
onto it.

The third 
example is the skein space of a disk with $2n$ points on the boundary.
If the disk is viewed as a rectangle with $n$ points on one side and $n$ on 
the opposite, then we can define a multiplication rule on the skein space by 
juxtaposing rectangles, obtaining the Temperley-Lieb algebra $TL_n$. Recall
that $TL_n$ is generated by the elements $1,e_1, e_2, \cdots, e_{n-1}$, where
$e_i$ is described in Fig. 3.2. 

\begin{figure}[htbp]
\centering
\leavevmode
\epsfxsize=1.7in
\epsfysize=1.2in
\epsfbox{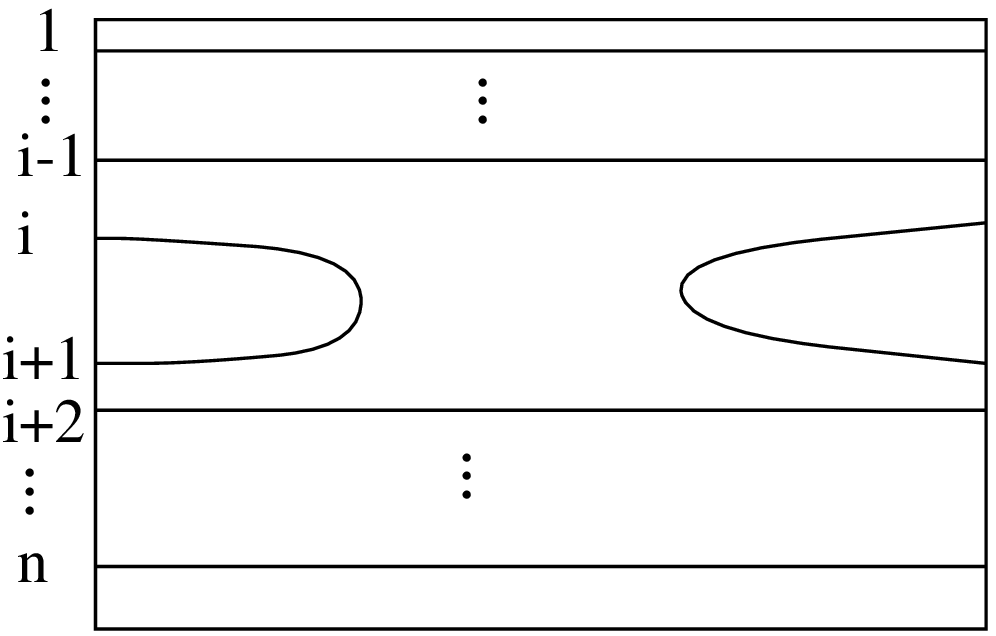}

Fig. 3.2.
\end{figure}

There exists a map from $TL_n$ to ${\cal S}({\bf R}^2)$ obtained by
closing the elements in $TL_n$ by $n$ parallel arcs. This map plays
the r{\^{o}}le of a quantum trace. It splits in a canonical way as
$TL_n\rightarrow {\cal S}(A)\rightarrow {\cal S}({\bf R}^2)$ by first 
closing the
elements in an annulus and then including them in a plane.

At this moment we recall the definition of the Jones-Wenzl idempotents [We]. 
They are of great importance for our construction, since they mimic
the behavior of the finite  dimensional irreducible representations
from the 
Reshetikhin-Turaev theory [RT]. For this let $r>1$ be an integer (which
will be the level of our TQFT). Let $A=e^{i\pi /(2r)}$. Recall that for each 
$n$ one denotes by $[n]$ the quantized integer $(A^{2n}-A^{-2n})/(A^2-A^{-2}).$

The Jones-Wenzl idempotents are the unique elements $f^{(n)}\in TL_n$,
$0\leq n\leq r-1$, that satisfy the following properties:

1) $f^{(n)}e_i=0=e_if^{(n)}$, for $0\leq i\leq n-1$,

2) $(f^{(n)}-1)$ belongs to the algebra generated by $e_1,e_2,\cdots ,e_{n-1}$,

3) $f^{(n)}$ is an idempotent,

4) $\Delta _n=(-1)^n[n+1]$

\noindent where $\Delta _n $ is the image of $f^{(n)}$
through the map $TL_n\rightarrow {\cal S}({\bf R}^2)$.

In the sequel we will have to work with the square root of  $\Delta  _n$ so we
make the notation $d_n=i^n\sqrt{[n+1]}$, thus $\Delta _n =d_n^2$.

Following [L1], in a diagram we will always denote $f^{(n)}$ by an empty coupon
(see Fig. 3.3).

\begin{figure}[htbp]
\centering
\leavevmode
\epsfxsize=1in
\epsfysize=0.6in
\epsfbox{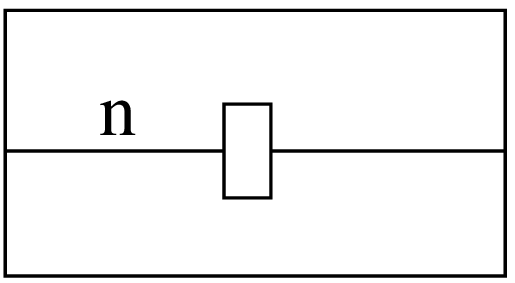}

Fig. 3.3.
\end{figure}

The image of $f^{(n)}$ 
through the map $TL_n\rightarrow {\cal S}(A)$ will be denoted
by $S_n(\alpha )$
We will also need the element $\omega \in {\cal S}(A)$,
$\omega =\sum_{n=0}^{r-2}d_n^2S_n(\alpha )$. Given a link diagram $L$ in the 
plane, whenever we label one of its components by $\omega $ we actually mean 
that we inserted $\omega $ in the way described in the definition of 
$<\cdot ,\cdot ,\cdots ,\cdot >_L$. 
Note that one can perform handle slides (also called second 
Kirby moves [Ki]), over components labeled by $\omega $ without changing
the value of the diagram (see [L1]).

Now we can define the basic data for a TQFT in level $r$, where $r$, as said,
is an integer greater than $1$.  Let ${\cal L}=\{0,1,\cdots ,r-2\}$.
Make the notation $X=(i\sqrt{2r})/(A^2-A^{-2})$, that is $X^2=\sum d_n^4=
<\omega >_U$, where $U$ is the unknot with zero framing (see [L1]).

Notice that by gluing two disks along the boundary we get a pairing map
${\cal S}({\bf D},2n)\times {\cal S}({\bf D},2n)\rightarrow {\cal S}(S^2)=
{\bf C}$,
hence we can view ${\cal S}({\bf D},2n)$ as a set of functionals acting on
the skein space of the exterior. In what follows, whenever we mention the 
skein space of a disk, we will always mean the skein space as a set of
functionals in this way. For example this will enable us
to get rid of the diagrams that have a strand labeled by $r-1$ (see
also [L1], [K2]). The point of view is similar to that of
factoring by the bad part of a representation (the one of quantum trace $0$) in
the Reshetikhin-Turaev setting.

To a disk with boundary labeled by $0$ we associate the vector space
$V_0$ which is the skein space of a disk with no points on the boundary.
Of course for any other label $a$ we put $V_a=0$. It is obvious that 
$V_0={\bf C}$. We let $\beta _0$ be the empty diagram.

\begin{figure}[htbp]
\centering
\leavevmode
\epsfxsize=2.9in
\epsfysize=1.1in
\epsfbox{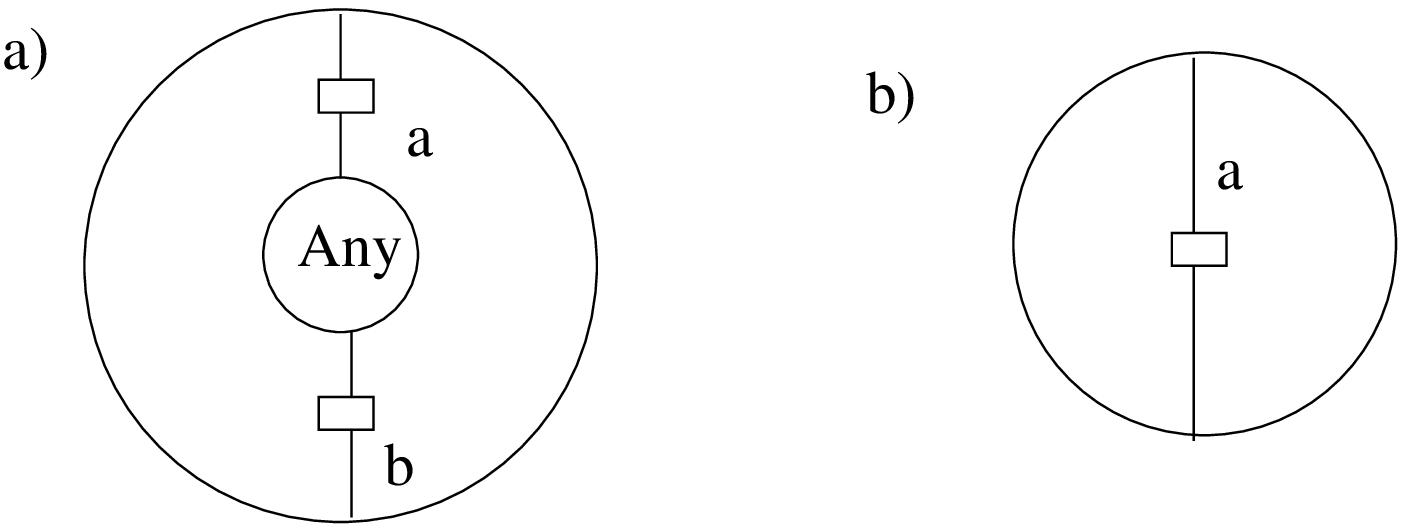}

Fig. 3.4.
\end{figure}

To an annulus with boundary components labeled by $a$ and $b$ we associate the 
vector space $V_{ab}$ which is the subspace of ${\cal S}({\bf D},a+b)$ 
spanned by 
all diagrams of the form indicated in Fig. 3.4. a), where in the smaller
disk can be inserted any diagram from 
${\cal S}({\bf D},a+b)$. The first 
condition in the definition of the Jones-Wenzl
idempotents implies that $V_{ab}=0$ if $a\neq b$ 
and $V_{aa}$ is one dimensional
and is spanned by the diagram from Fig. 3.4. b). We will denote by
$\beta _{aa}$ this diagram multiplied by $1/d_a$, where we recall that
$d_a=i^a\sqrt{[a+1]}$. The element $\beta _{aa}$ has the property that paired 
with itself on the outside gives 1.

To a pair of pants with boundary components labeled by $a$, $b$, and $c$
we put into correspondence the space $V_{abc}$, which is the space
spanned by all diagrams of the form described in Fig. 3.5. a), where in the
inside disk we allow any diagram from ${\cal S}({\bf D},a+b+c)$.

\begin{figure}[htbp]
\centering
\leavevmode
\epsfxsize=4.7in
\epsfysize=1.3in
\epsfbox{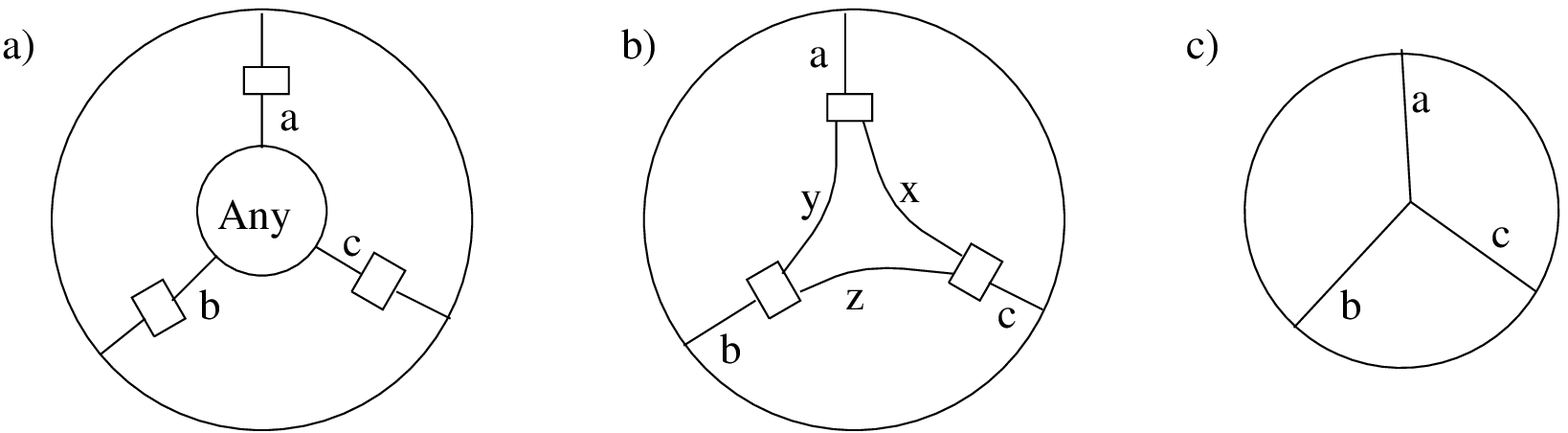}

Fig. 3.5.
\end{figure}

The reader will notice that there is some
ambiguity in this definition. To make it rigorous, we have to mark a point 
on the circle, from which all points are counted. We will keep this in
mind although we will no longer mention it.

The results from [K2] and [L1] show that $V_{abc}$ can either be one
dimensional or it is equal to zero. The triple $(a,b,c)$ is said to be
admissible if $V_{abc}\neq 0$. This is exactly the case when $a+b+c$ is
even, $a+b+c\leq 2(r-2)$ and $a\leq b+c$, $b\leq a+c$, $c\leq a+b$.
In this case the space $V_{abc}$ is spanned by the triad introduced by
Kauffman [K2] which is described in Fig. 3.5. b). Here the numbers 
$x$, $y$, $z$ satisfy $a=x+y$, $b=y+z$, $c=z+x$.

\begin{figure}[htbp]
\centering
\leavevmode
\epsfxsize=5in
\epsfysize=2in
\epsfbox{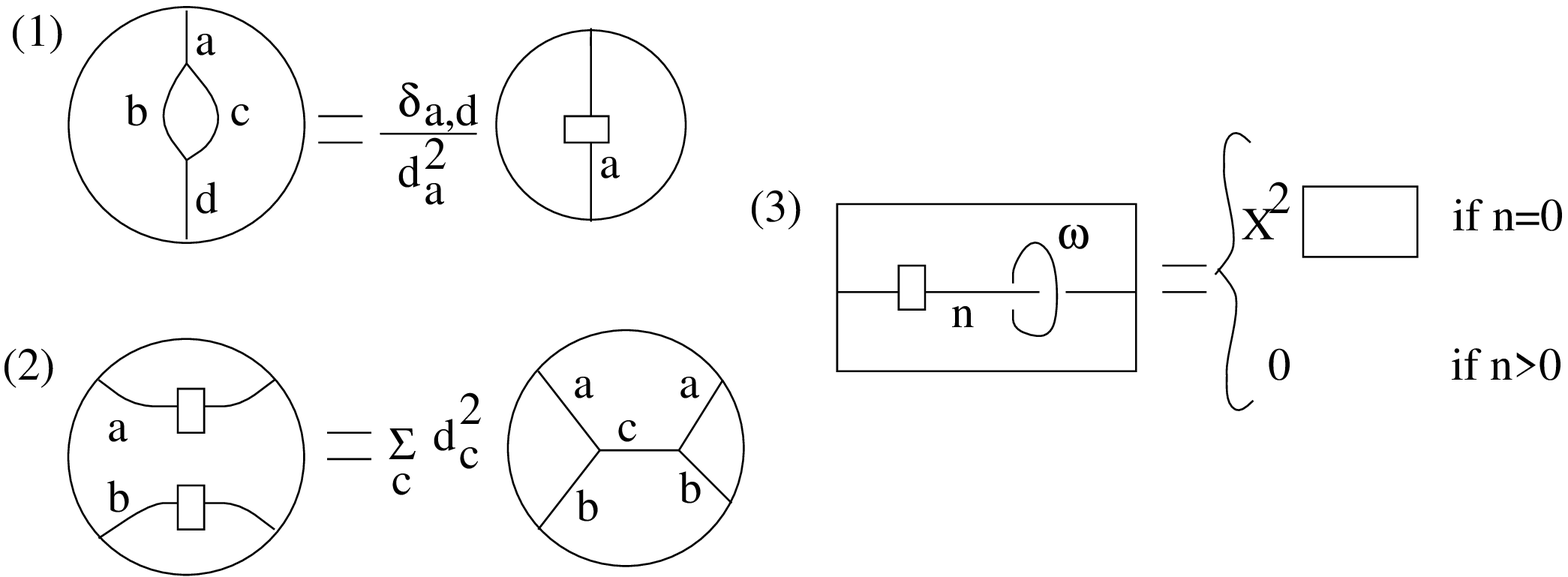}

Fig. 3.6.
\end{figure}

In [L2] it is shown that if we pair the diagram from Fig. 3.5. b) with the
one corresponding to $V_{acb}$ on the outside we get the complex number 
$\theta (x+y,y+z,z+x)=(\Delta _{x+y+z}!\Delta _{x-1}!\Delta _{y-1}!\Delta
_{z-1}!)/$ $(\Delta _{y+z-1}!\Delta _{z+x-1}!\Delta _{x+y-1}!)$, where
$\Delta _n=\Delta _1\Delta_2\cdots \Delta _n$ and $\Delta _{-1}=1$.
Thus if we denote by $\beta _{abc}$ the product of this diagram 
with $(d_{x+y+z}!d_{x-1}!d_{y-1}!d_{z_1}!)/(d_{y+z-1}!d_{z+x-1}!d_{x+y-1}!)=
1/\sqrt{\theta (a,b,c)}$ (with the same convention for factorials),
then $\beta _{abc}$ paired on the outside with $\beta _{acb}$ will give $1$.

In diagrams, whenever we have a 
$\beta _{abc}$ we make the notation from Fig. 3.5. c). This 
notation is different from the one  with a dot in the
middle from [L1] , in the sense that we have a different normalization!
We prefer this notation because it will simplify diagrams in the future, so
whenever in a diagram we have a trivalent vertex, we consider that we
have a $\beta $ inserted there. In particular, a diagram that looks like the 
Greek letter $\theta $ will be equal to $1$ in ${\cal S}({\bf R}^2)$.
The elements $\beta _{abc}$ are the analogues of the quantum Clebsch-Gordan
coefficients.

In the sequel we will need the three identities described in Fig. 3.6, whose
proofs can be found in [L1]. Here $\delta _{ad}$ is the Kronecker symbol.

\begin{figure}[htbp]
\centering
\leavevmode
\epsfxsize=2.4in
\epsfysize=2.1in
\epsfbox{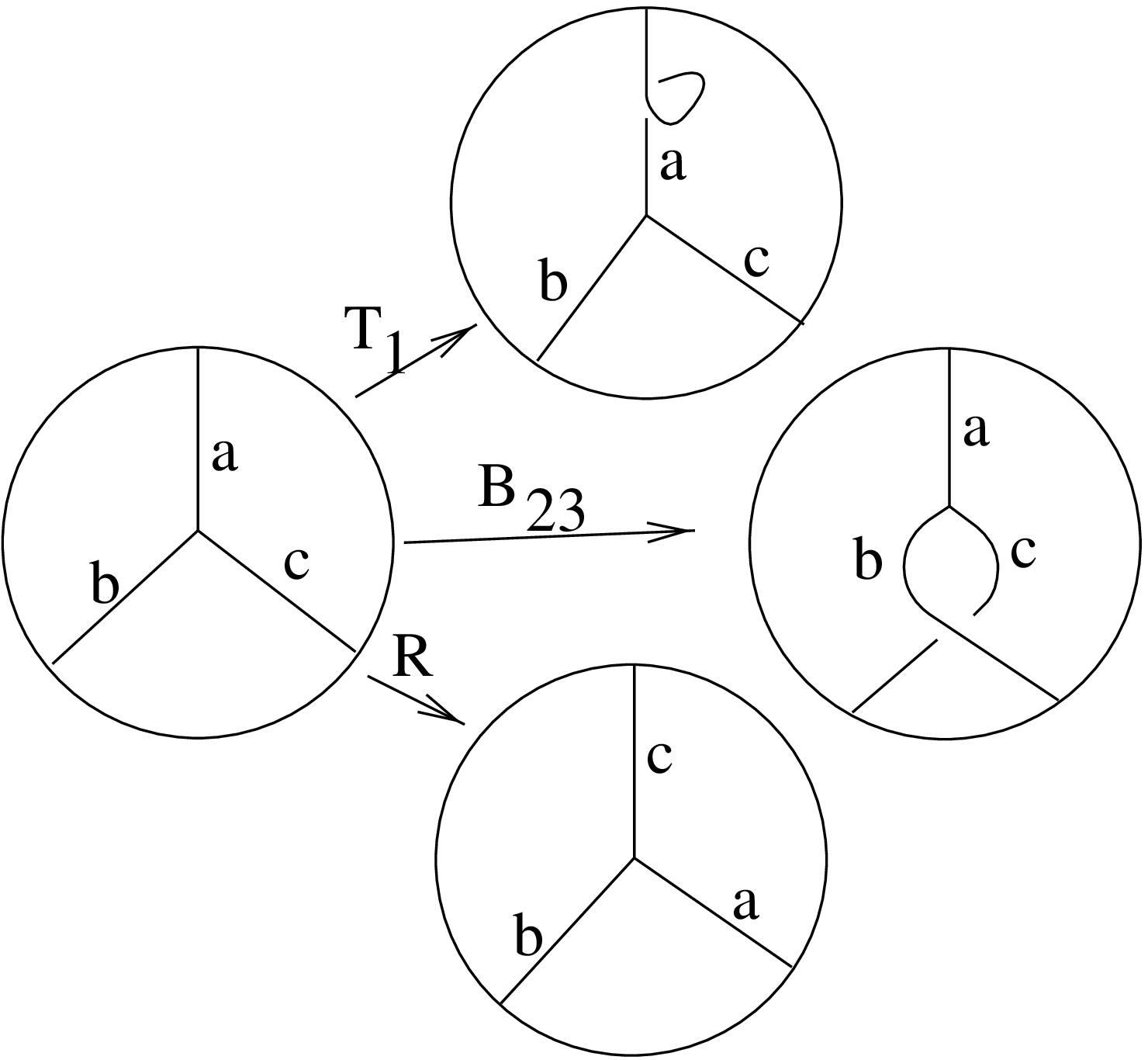}

Fig. 3.7.
\end{figure}

Let us define the dual spaces. It is natural to let the
dual of $V_0$ to be $V_0$, that of $V_{aa}$ to be $V_{aa}$, and that
of $V_{abc}$ to be $V_{acb}$. However the pairings will look peculiar.
This is due to the fact that
 we want the mapping cylinder axiom to be satisfied.
So we let $<,>:V_0\times V_0\rightarrow {\bf C}$ be defined by
$<\beta _0,\beta _0>=1$, $<,>:V_{aa}\times V_{aa}
\rightarrow {\bf C}$ be defined by
$<\beta _{aa},\beta _{aa}>=X/d_a^2$, 
and $<,>:V_{abc}\times V_{acb}\rightarrow {\bf C}$ be defined by
$<\beta _{abc},\beta _{acb}>=X^2/(d_ad_bd_c)$.

\begin{figure}[htbp]
\centering
\leavevmode
\epsfxsize=4in
\epsfysize=0.9in
\epsfbox{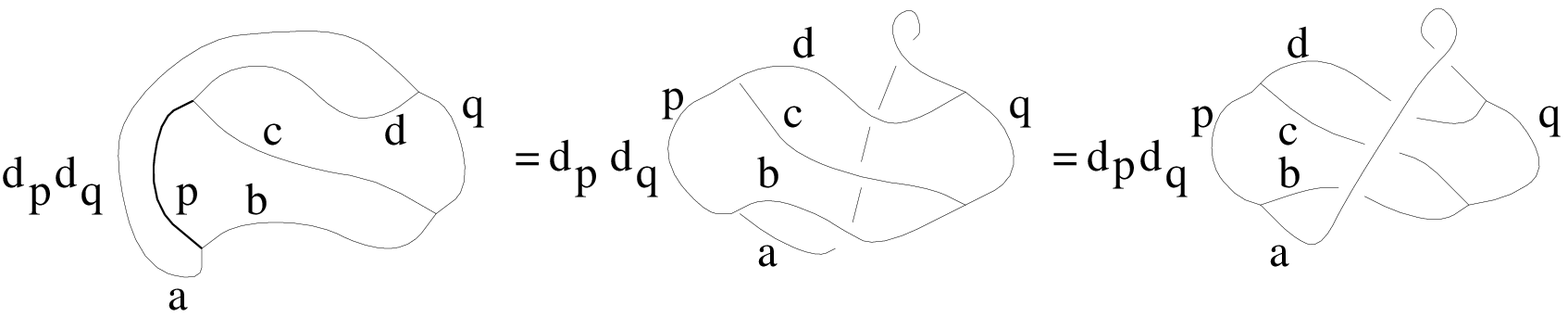}

Fig. 3.8.
\end{figure}

Before we define the morphisms associated to the elementary moves we make
the convention that for any e-morphism $f$ we will denote $V(f)$ also by 
$f$.

The morphisms corresponding to the three elementary moves on a pair of pants 
are described in Fig. 3.7.
Further, we let $F:\bigoplus _pV_{pab}\otimes V_{pcd}\rightarrow 
\bigoplus _pV_{qda}\otimes V_{qbc}$ be defined by 
$
F\beta _{pab}\otimes \beta _{pcd}=\sum_q f_{abcdpq}
\beta _{qda}\otimes \beta _{qbc}
$
the coefficients $f_{abcdpq}$ being given by any of the three equal diagrams 
from Fig. 3.8. Note that $f_{abcdpq}=d_p^{-1}d_q\{^{bcp}_{adq}\}$ 
where $\{^{bcp}_{adq}\}$ are the  6j-symbols.

\begin{figure}[htbp]
\centering
\leavevmode
\epsfxsize=2.7in
\epsfysize=0.9in
\epsfbox{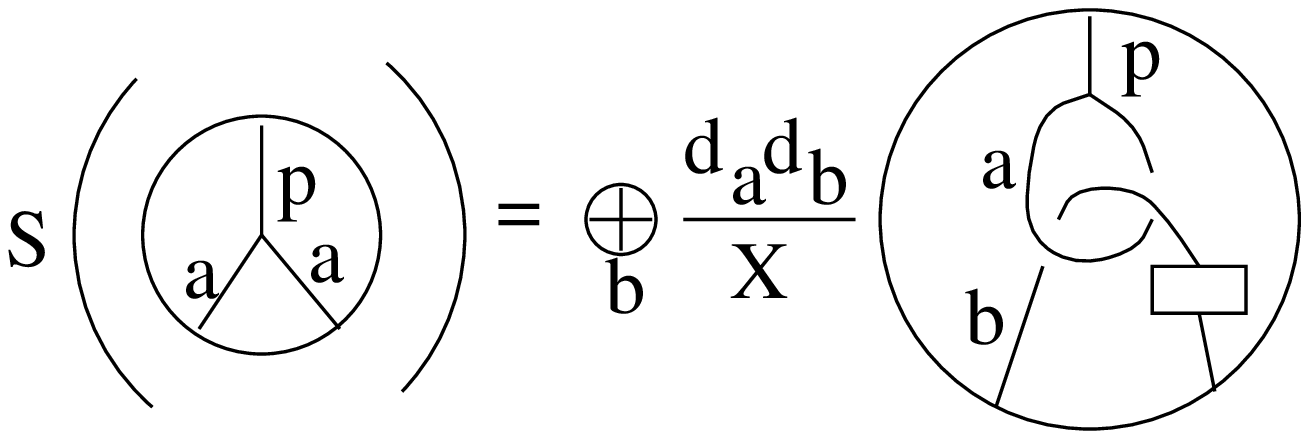}

Fig. 3.9.
\end{figure}

Also the map $S:\bigoplus V_{paa}\rightarrow \bigoplus V_{pbb}$ is described in
Fig. 3.9.

\begin{figure}[htbp]
\centering
\leavevmode
\epsfxsize=2in
\epsfysize=0.5in
\epsfbox{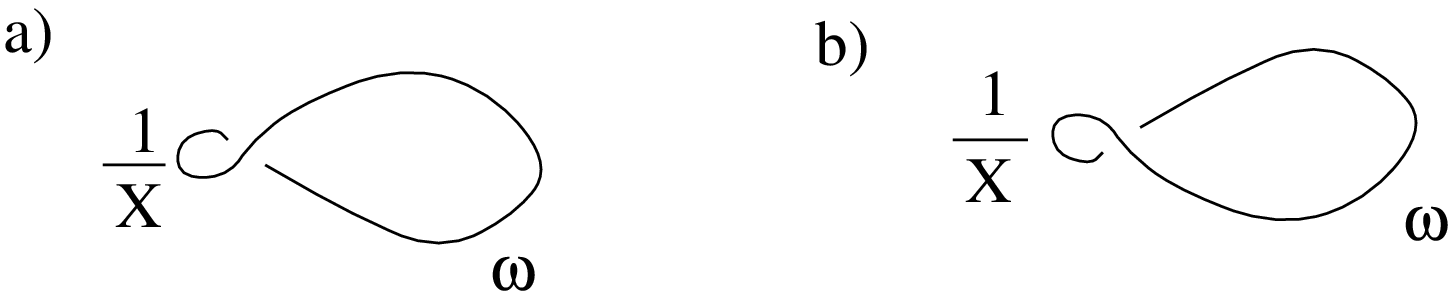}

Fig. 3.10.
\end{figure}

The maps $A$, $D$ and $P$ are given by relations of the form $A(x\otimes \beta 
_{aa})=x$, $D(\beta _{aa0}\otimes \beta _0)=\beta _{aa}$ and $P(x\otimes y)=y
\otimes x$. The map $C$ is the multiplication by the value of the diagram 
described in Fig. 3.10. a). Note that Lemma 4 in [L1] implies that the
inverse of $C$ is the multiplication by the diagram from Fig. 3.11. b).
Finally, $S(a)=d_a^2/X$, $a\in {\cal L}$.

\medskip

\ul{Remark} The reader should note that the crossings from all these diagrams 
are negative.  We make this choice because, returning to the analogy
with vector spaces, all the maps we defined behave like changes of basis
rather than like morphisms.

\bigskip

\begin{center}
{\bf 4. THE COMPATIBILITY CONDITIONS}
\end{center}

\bigskip

In order for the basic data to give rise to a well defined TQFT, it has to
satisfy certain conditions. A list of such conditions has been exhibited
in [Wa], by making use of techniques of Cerf theory similar
to those from [HT]. The first group of relations, the so called 
Moore-Seiberg equations, are the conditions that have to be satisfied in order 
for the modular functor to exist. They are as follows:

\medskip

1. at the level of a pair of pants:

\ a) $T_1B_{23}=B_{23}T_1$,\ $T_2B_{23}=B_{23}T_3$, \ $T_3B_{23}=B_{23}T_2$,
where $T_2=RT_1R^{-1}$ and $T_3=R^{-1}T_1R$,

\ b) $B_{23}^2=T_1T_2^{-1}T_3^{-1}$,

\ c) $R^3=1$,

\ d) $RB_{23}R^2B_{23}RB_{23}R^2=1$,

\begin{figure}[htbp]
\centering
\leavevmode
\epsfxsize=3.2in
\epsfysize=0.8in
\epsfbox{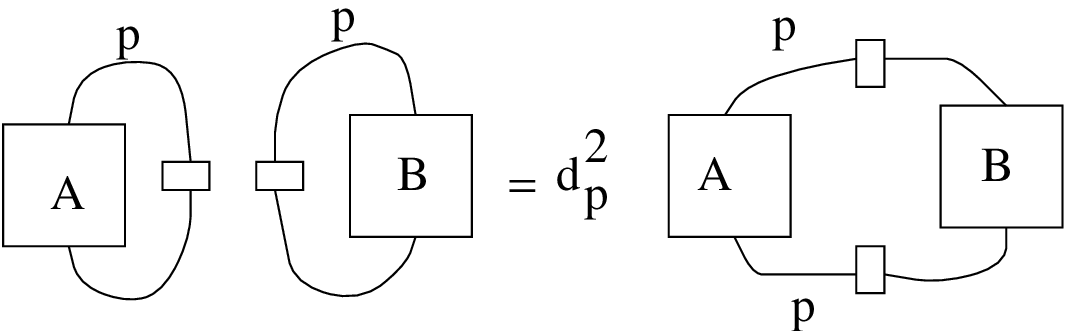}

Fig. 4.1.
\end{figure}

2. relations defining  inverses:

\ a) $P^{(12)}F^2=1$,

\ b) $T_3^{-1}B_{23}^{-1}S^2=1$,

3. relations coming from ``codimension 2 singularities'':

\ a) $P^{(13)}R^{(2)}F^{(12)}R^{(2)}F^{(23)}R^{(2)}
F^{(12)}R^{(2)}F^{(23)}R^{(2)}F^{(12)}=1$,

\ b) $T_3^{(1)}FB_{23}^{(1)}FB_{23}^{(1)}FB_{23}^{(1)}=1$,

\ c) $C^{-1}B_{23}^{-1}T_3^{-2}ST_3^{-1}ST_3^{-1}S=1$,

\ d) $R^{(1)}(R^{(2)})^{-1}FS^{(1)}FB_{23}^{(2)}
B_{23}^{(1)}=FS^{(2)} T_3^{(2)}(T_1^{(2)})^{-1}B_{23}^{(2)}
F$,

4. relations involving annuli and disks:

\ a) $F(\beta _p^{mn}\otimes 
\beta _p^{p0})=\beta _m^{0m}\otimes \beta _m^{np}$,

\ b) $A^{(12)}_2D_3^{(13)}=
D_2D_3^{(13)}$,

\ c) $A^{(12)}A^{(23)}=
A^{(23)}A^{(12)}$,

\begin{figure}[htbp]
\centering
\leavevmode
\epsfxsize=5in
\epsfysize=1.7in
\epsfbox{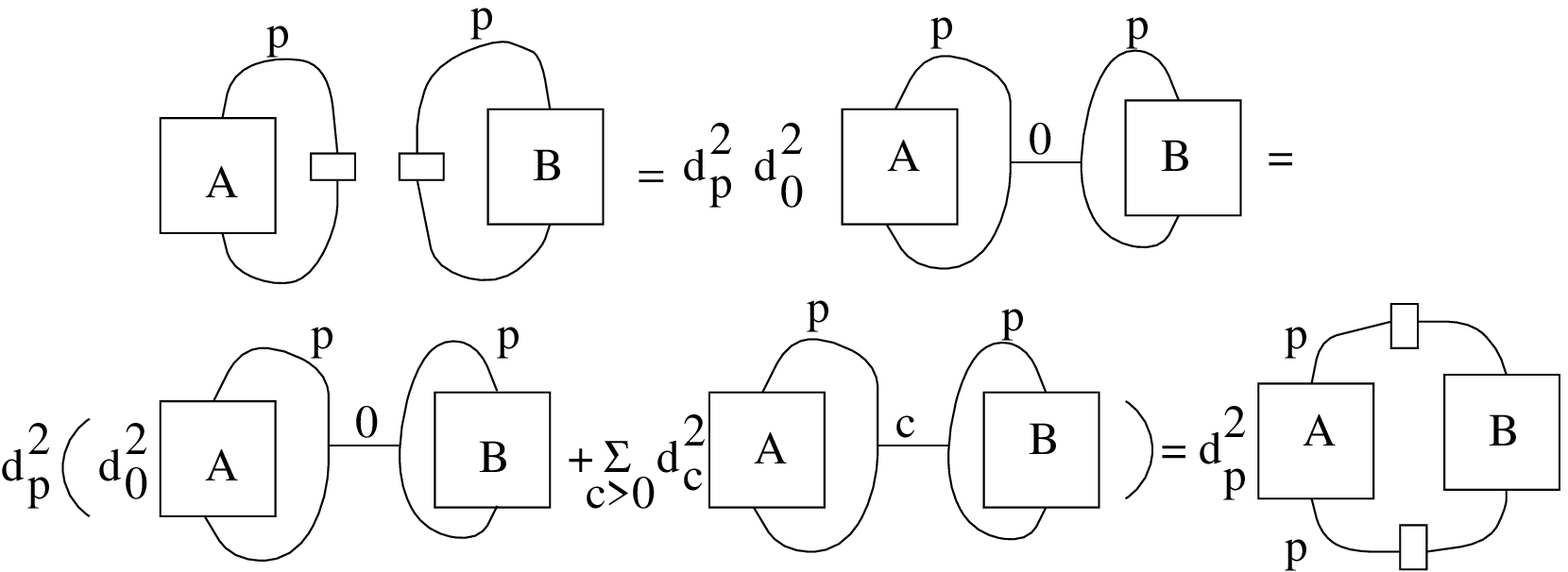}

Fig. 4.2.
\end{figure}

5. relations coming from duality:

\ -for any elementary move $f$, one must have $f^+=\bar{f}$, where
$f^+$ is the adjoint of $f$with respect to the pairing, and
$\bar{f}$ is the morphism induced by $f$ on the surface with reversed 
orientation,

6. relations expressing the compatibility between the pairing, and moves
$A$ and $D$:

\ a) $<\beta _m^m,\beta _m^m>=S(m)^{-1}$

\ b) $<\beta _m^{m0},\beta _m^{m0}>=S(0)^{-1}S(m)^{-1}$.

\medskip

In addition one also has to consider two conditions that guarantee that the
partition function is well defined.

7. a) $S(m)=S_{0m}$
where $[S_{xy}]_{x,y}$  is the matrix of move $S$ on the torus (which can be
thought as the punctured torus capped with a disk),

\ b) $F(\beta _0^{mm}
\otimes\beta _0^{nn})=\bigoplus S(m)^{-1}S(n)^{-1}id_{pmn}$
where $id_{pmn}$ is the identity matrix in $(V_{pmn})^*\bigotimes V_{pmn}$.

\medskip

In all these relations, the superscripts in parenthesis indicate the index
of the elementary surface(s) on which the map acts, and the subscripts
indicate the number of the boundary component.

\begin{figure}[htbp]
\centering
\leavevmode
\epsfxsize=5in
\epsfysize=0.7in
\epsfbox{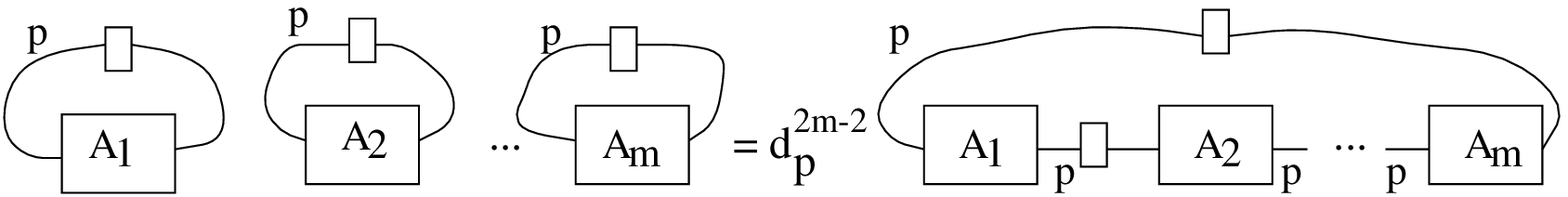}

Fig. 4.3.
\end{figure}

\begin{figure}[htbp]
\centering
\leavevmode
\epsfxsize=5.8in
\epsfysize=1in
\epsfbox{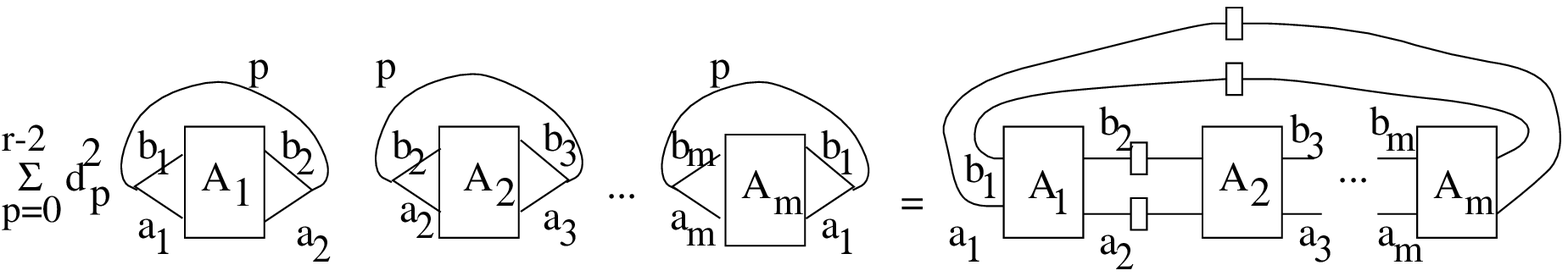}

Fig. 4.4.
\end{figure}

We will prove that our basic data satisfies these relations. For the
proof we will need a contraction formula similar to the tensor
contraction formula that one encounters in the case of TQFT's based on 
representations of Hopf algebras (see [T], [FK], [G1]).

\medskip

\ul{Lemma 4.1.} For any $A,B\in TL_p$ the equality from Fig. 4.1 holds.

\medskip

\ul{Proof:} The proof is contained in Fig. 4.2. In this chain of equalities
the first one follows from the way we defined the $\beta$'s,
 the second one holds because the sum that appears
in the third term is zero (by the first property of Jones-Wenzl
idempotents, since such an idempotent lies on the strand labeled by $c$;
more explanations about this phenomenon can be found in [L1] and [R]), and 
the last equality follows from identity (2) in Fig. 3.6.$\Box$

\medskip

\ul{Lemma 4.2.} If $A_1,A_2,\cdots ,A_m\in TL_p$ then the identity from 
Fig. 4.3 holds.

\medskip

\ul{Proof:} Follows by induction from Lemma 4.1.$\Box$.

\medskip

\ul{Theorem 4.1.} Suppose that $A_i\in {\cal S}({\bf D}, a_i+b_i+a_{i+1}+
b_{i+1})$, $i=1,2,\cdots, m$, where $a_i$ and $b_i$ are integers with 
$a_{m+1}=a_1$ and $b_{m+1}=b_1$. Then the identity described in Fig. 4.4
holds.

\medskip

\ul{Proof:} By Lemma 4.2, the left hand side is equal to the expression 
described in Fig. 4.5. a).

\begin{figure}[htbp]
\centering
\leavevmode
\epsfxsize=6.3in
\epsfysize=1in
\epsfbox{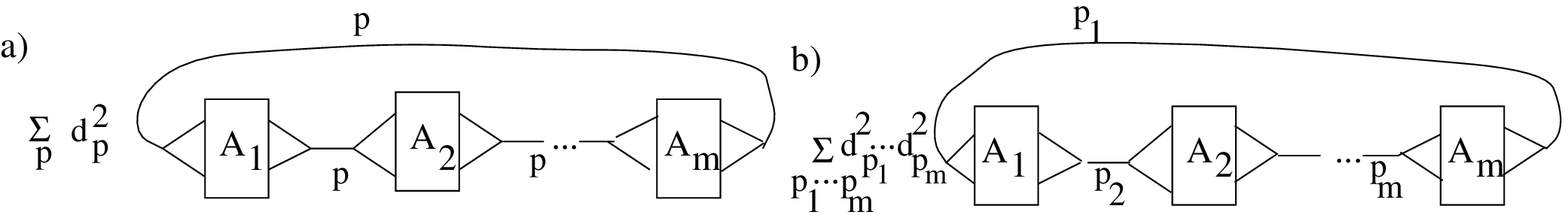}

Fig. 4.5.
\end{figure}

On the other hand, if $p\neq q$, by using the identity (2) from Fig. 3.6,
we get the chain of equalities from Fig. 4.6, where the last one follows 
from the fact that on the strand labeled by $c$ there is a Jones-Wenzl
idempotent and using the first property of these idempotents.

\begin{figure}[htbp]
\centering
\leavevmode
\epsfxsize=3.1in
\epsfysize=.8in
\epsfbox{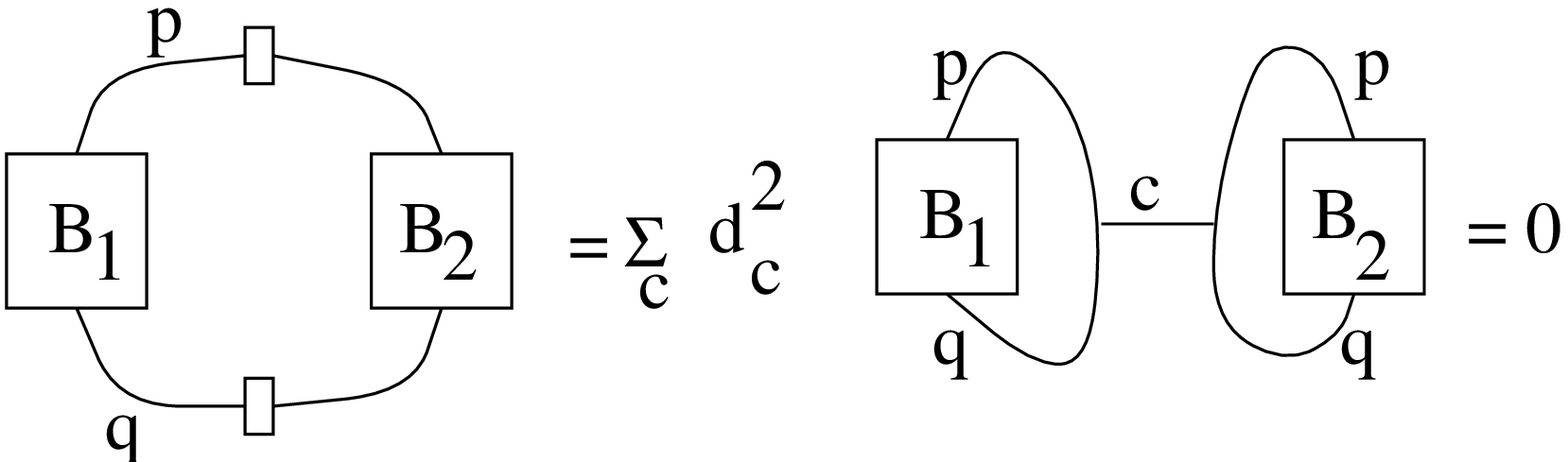}

Fig. 4.6.
\end{figure}

As a consequence of this fact we get that our expression is equal to the one
from Fig. 4.5. b), and then by applying the identity (2) from Fig. 3.6 several
times we get the desired  result.$\Box$

\medskip

\begin{figure}[htbp]
\centering
\leavevmode
\epsfxsize=3.8in
\epsfysize=0.8in
\epsfbox{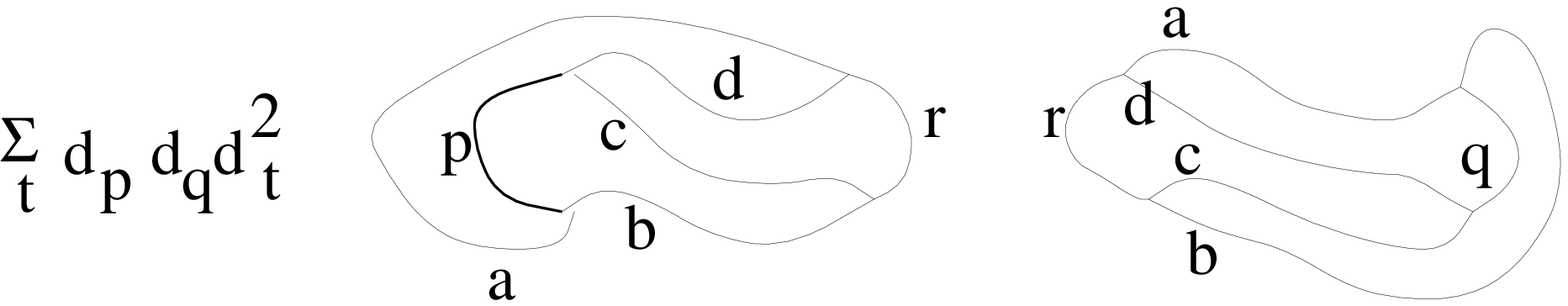}

Fig. 4.7.
\end{figure}

We can proceed with proving the compatibility conditions. The proofs
are similar to the ones in [FK] and [G1], but one should note that
they are simpler. First, the relations on a pair of pants are obviously
satisfied. This can be seen at first glance for 1.a) and 1.c), then
1.d) is the third Reidemeister move, and 1.b) is equivalent to 1.c) (see
[FK] or Chap. VI  in [T]).

For the proof of 2.a) we write $FPF\beta _{pab}\otimes \beta _{pcd}=
\sum_q  c_{abcdpq}\beta _{qab}\otimes \beta _{qcd}$. Since we have
a matrix multiplication here we see that the coefficient $c_{abcdpq}$
is given by the diagram from Fig. 4.7.

By using Theorem 4.1 this becomes the expression from Fig. 4.8. Using 
identity (1) from Fig. 3.6. wee see that this is equal to $\delta _{pq}$
multiplied by the Greek letter $\theta $ diagram, therefore is equal to 
$\delta _{pq}$ and the identity is proved.

\begin{figure}[htbp]
\centering
\leavevmode
\epsfxsize=1.9in
\epsfysize=1.1in
\epsfbox{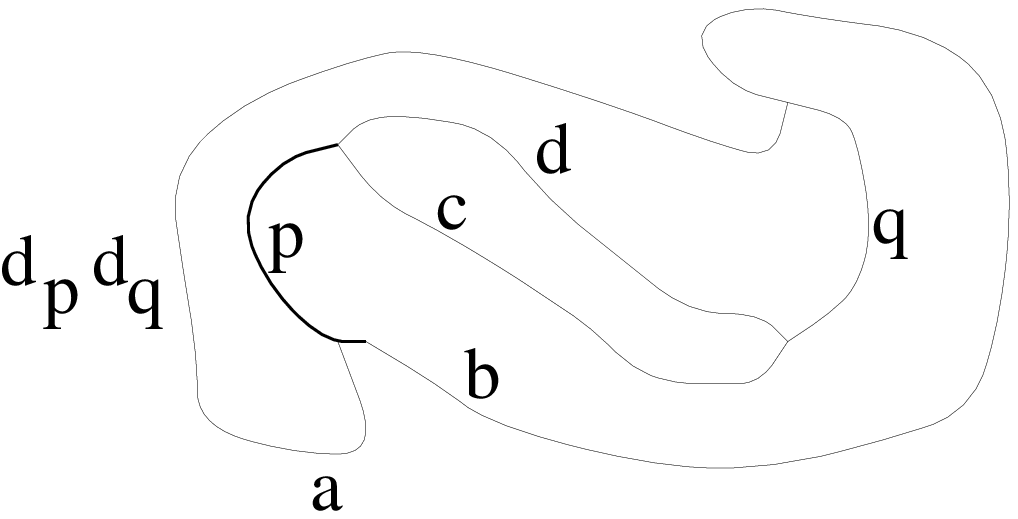}

Fig. 4.8.
\end{figure}

For 2.b) we have that $T_3^{-1}B_{23}^{-1}S^2\beta _{paa}$ is equal to 
the first term in Fig. 4.9. We get the chain of equalities from this
figure by pulling first
the strand labeled by  $\omega $
 down and using  the identity (2) from Fig. 3.6, and then using identity
(3) from Fig. 3.6. The last term is equal to $\beta _{aa}$.

\begin{figure}[htbp]
\centering
\leavevmode
\epsfxsize=4.8in
\epsfysize=1.1in
\epsfbox{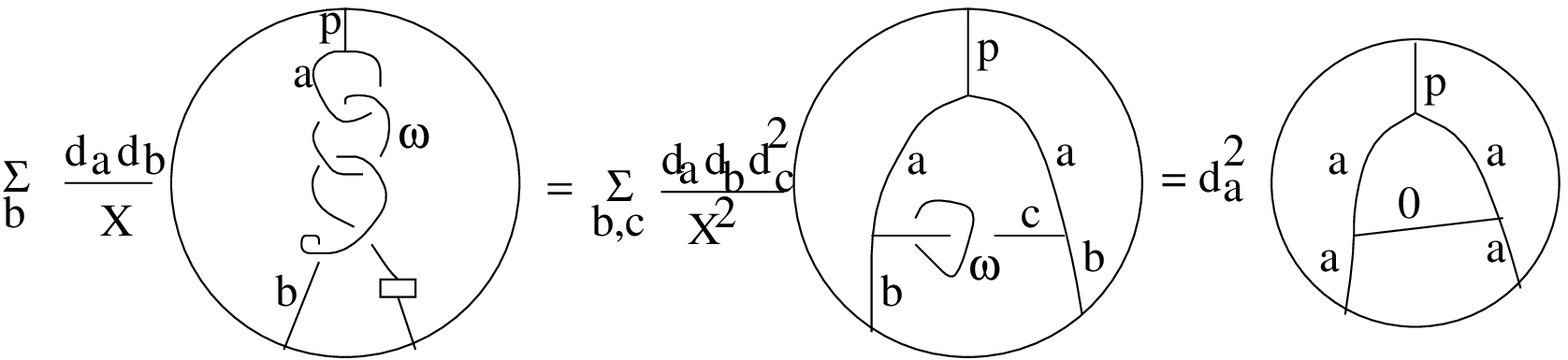}

Fig. 4.9.
\end{figure}

\begin{figure}[htbp]
\centering
\leavevmode
\epsfxsize=4.8in
\epsfysize=1.7in
\epsfbox{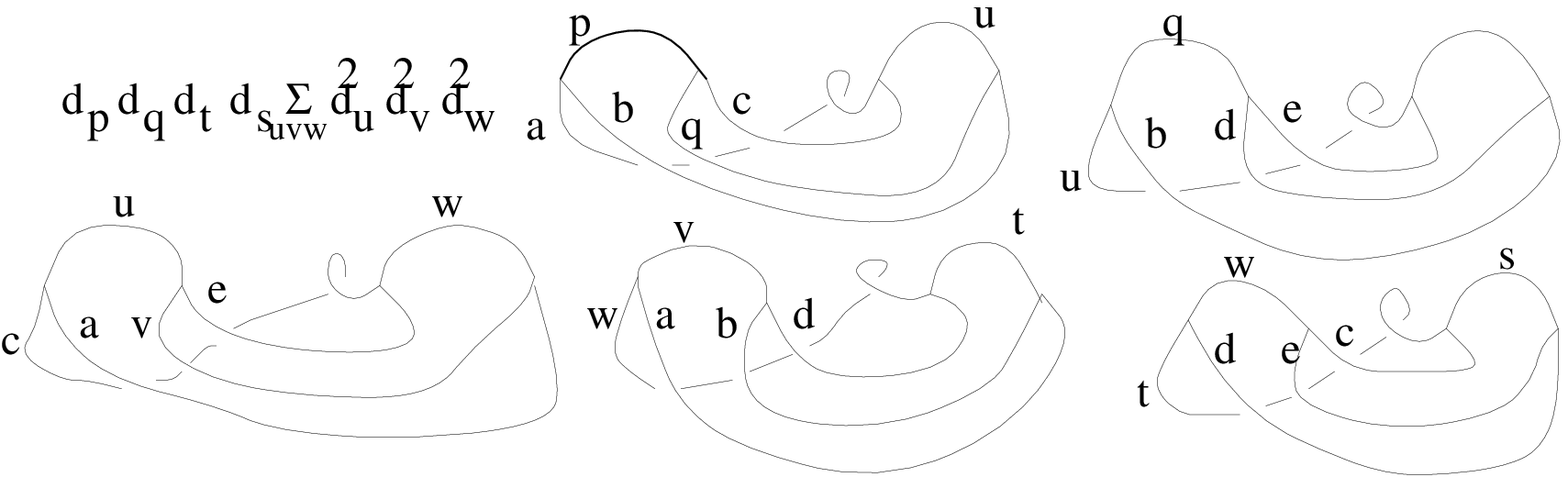}

Fig. 4.10.
\end{figure}

Now we describe the proof of the pentagon identity. We are interested in 
computing the coefficient of $\beta _{sde}\otimes \beta _{sct} \otimes
\beta _{rab}$ in $F^{(12)}R^{(2)}F^{(23)}R^{(2)}
F^{(12)}R^{(2)}F^{(23)}R^{(2)}F^{(12)}\beta _{pab}\otimes \beta _{pqc}
\otimes \beta _{qde}$. Again, by using the formula for matrix multiplication
we get that this coefficient is described in Fig. 4.10.

\begin{figure}[htbp]
\centering
\leavevmode
\epsfxsize=4.8in
\epsfysize=3.7in
\epsfbox{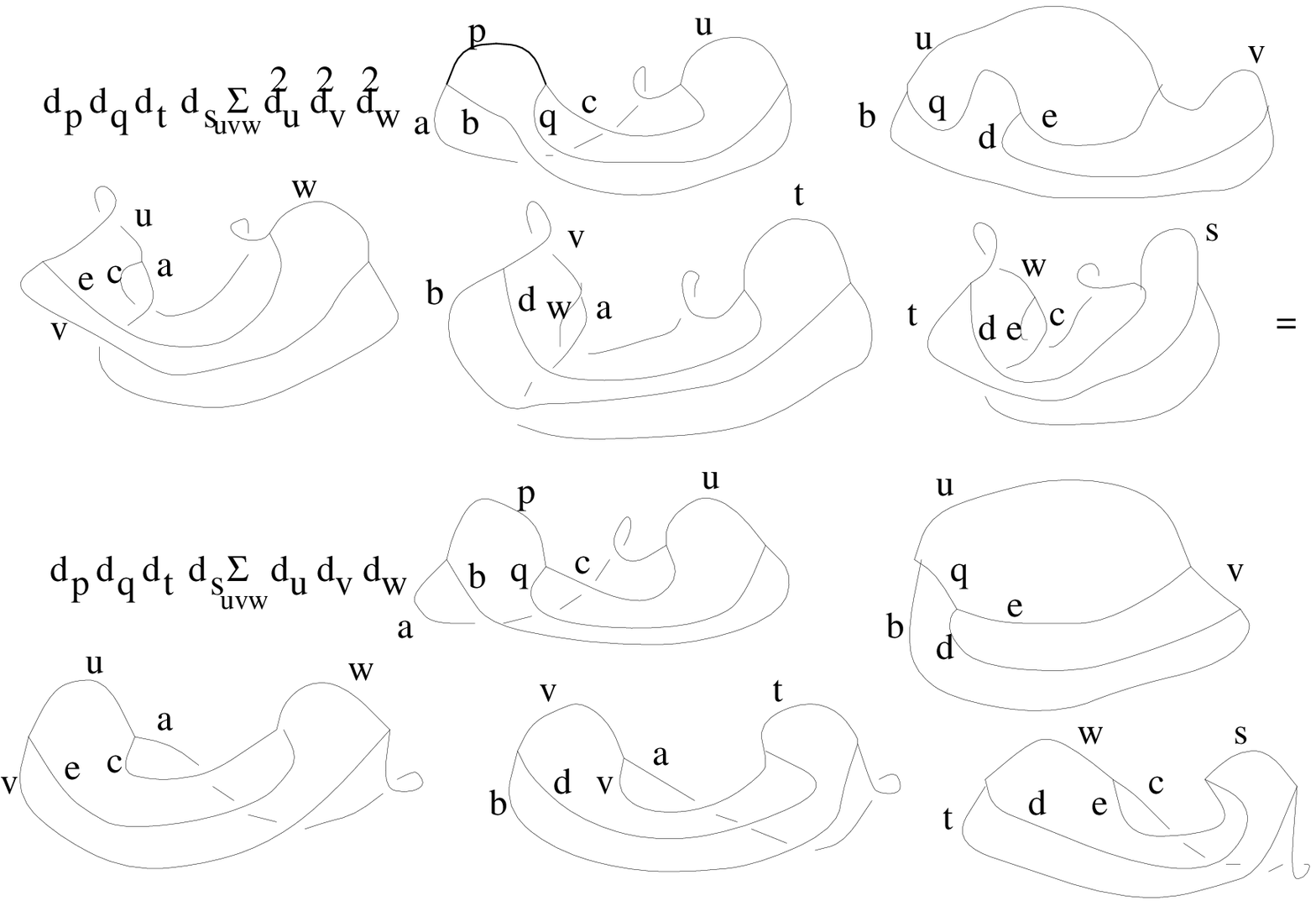}

Fig. 4.11.
\end{figure}

\begin{figure}[htbp]
\centering
\leavevmode
\epsfxsize=5.1in
\epsfysize=3.7in
\epsfbox{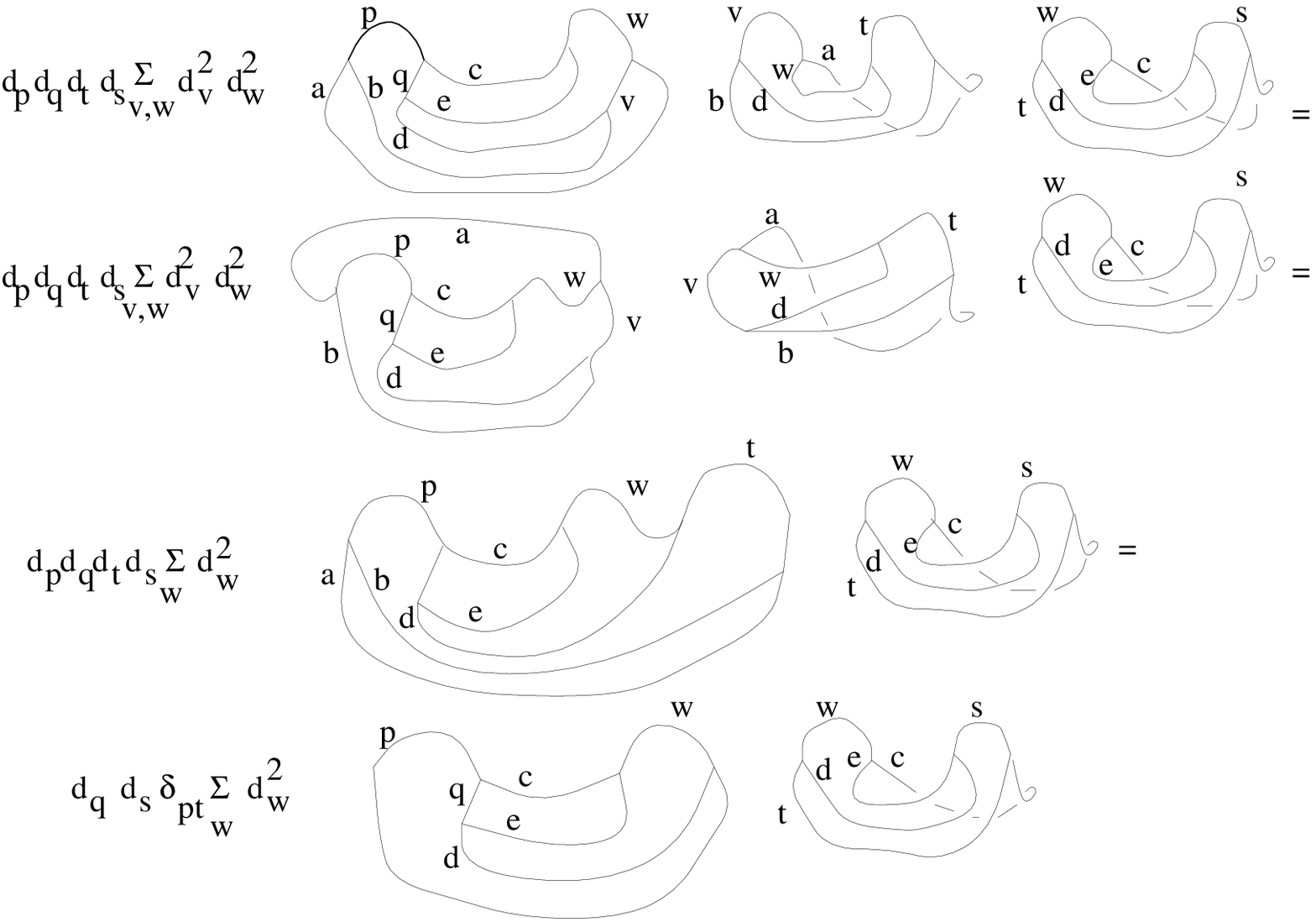}

Fig. 4.12.
\end{figure}

By doing a flip in the third, fourth
 and fifth factor we get the first term from
the equality shown in Fig. 4.11, which is further transformed into the 
second by applying three times 1.b).
Apply Theorem 
4.1 to contract with respect to
 $u$, then continue like in Fig. 4.12, namely 
pull the strand of $a$ over, then apply Theorem 4.1 for the sum over $v$
and then use for the last equality formula (1) in Fig. 3.6.
Finally, if we use Theorem 4.1 once more and then formula (1) in Fig. 3.6,
we get $\delta _{pt}\delta _{qs}$ times a diagram of the form of 
letter $\theta$. Hence the final answer is $\delta _{pt}\delta _{qs}$
and the identity is proved.

In order for the F-triangle to hold we have to show that the coefficient
of $\beta _{qab}\otimes \beta _{qcd}$ in 
$T_3^{(1)}FB_{23}^{(1)}FB_{23}^{(1)}FB_{23}^{(1)}\beta _{pab}\beta _{pcd}$
is $\delta _{pq}$. The coefficient is given in Fig. 4.13.

\begin{figure}[htbp]
\centering
\leavevmode
\epsfxsize=4.4in
\epsfysize=0.7in
\epsfbox{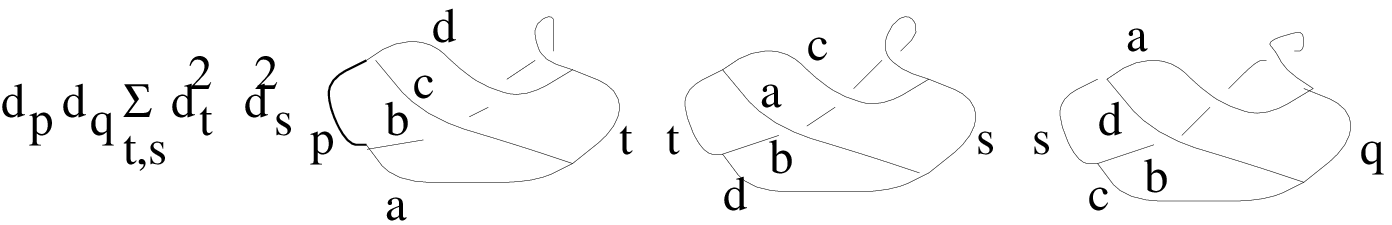}

Fig. 4.13.
\end{figure}

We transform the second factor as shown in Fig. 4.14 by first doing two 
flips and then using 1.b) twice. Then contract the product via
Theorem 4.1 to get the first term from the equality from Fig. 4.15, then
transform it into the second by using again 1.b).
As before, this is equal to $\delta _{pq}$.

\begin{figure}[htbp]
\centering
\leavevmode
\epsfxsize=4.2in
\epsfysize=0.8in
\epsfbox{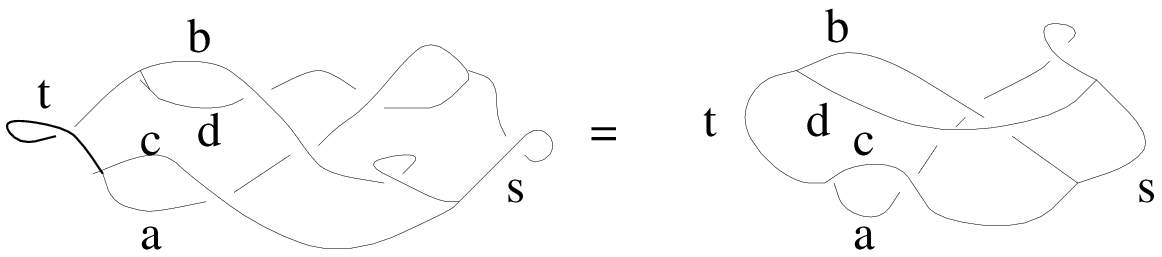}

Fig. 4.14.

\centering
\leavevmode
\epsfxsize=4in
\epsfysize=0.8in
\epsfbox{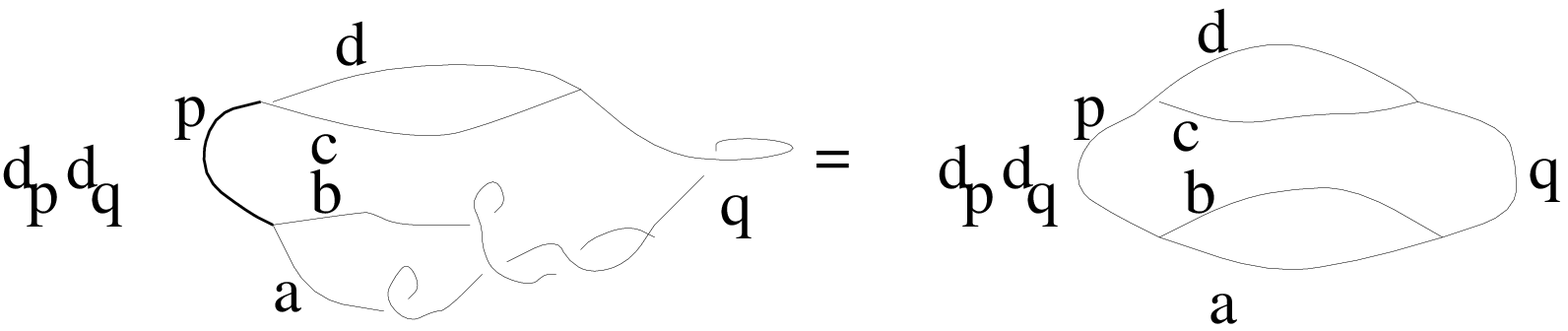}

Fig. 4.15.
\end{figure}

In the case of the S-triangle, it is not hard to see that 
$C^{-1}B_{23}^{-1}T_3^{-2}ST_3^{-1}ST_3^{-1}S\beta _{aa}$ 
is equal to the expression from Fig. 4.16. 
Lemma 3 in [L1] enables us to do Kirby moves over components labeled
by $\omega$, so we get the first term from Fig. 4.17, which is
equal to the second one by Lemma 4 in [L1]. From here we continue like
in the proof of 2.b).

\begin{figure}[htbp]
\centering
\leavevmode
\epsfxsize=2in
\epsfysize=1.2in
\epsfbox{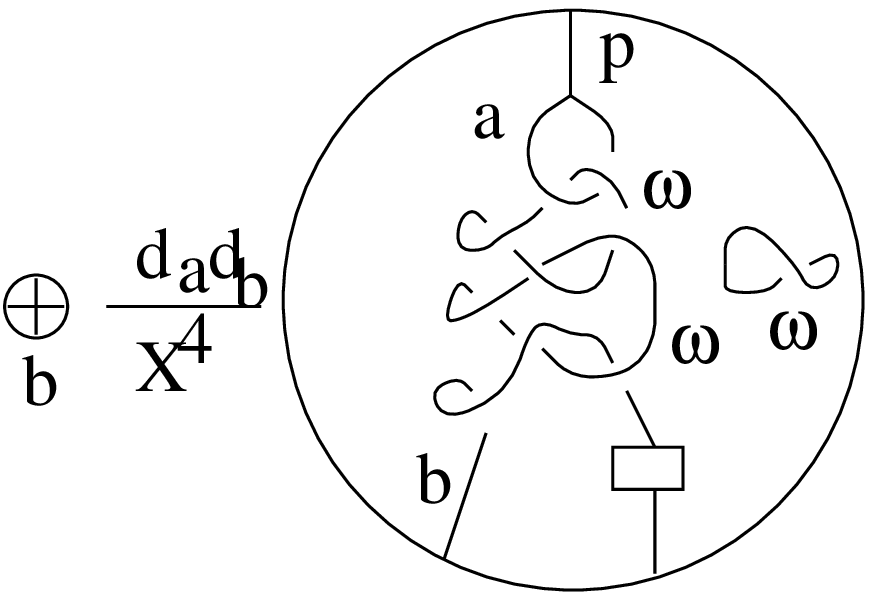}

Fig. 4.16.
\end{figure}

\begin{figure}[htbp]
\centering
\leavevmode
\epsfxsize=4in
\epsfysize=1.1in
\epsfbox{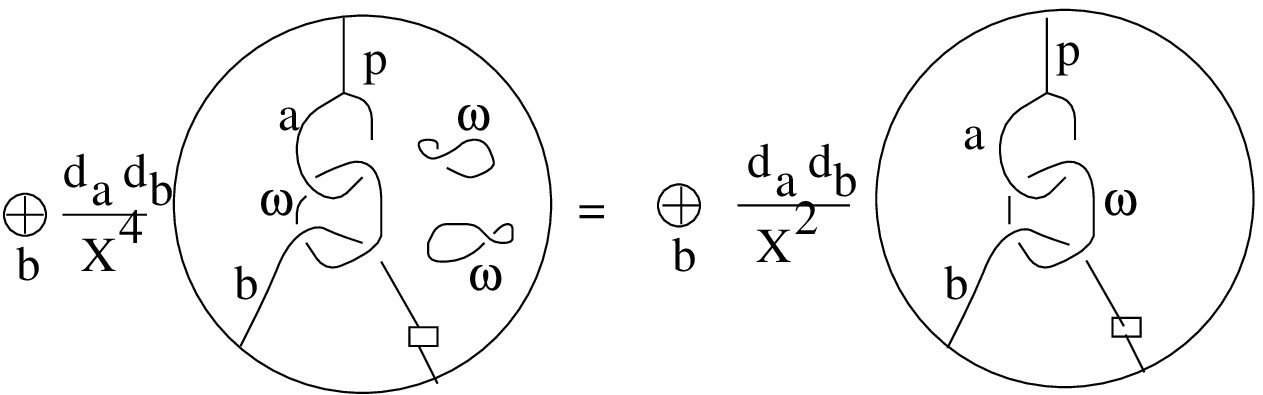}

Fig. 4.17.
\end{figure}

Let us prove 3.d). We have to show that the coefficient of
$\beta _{qdc}\otimes \beta _{qda}$ in
$FS^{(2)} T_3^{(2)}(T_1^{(2)})^{-1}B_{23}^{(2)}
F$ $\beta _{pab}\otimes \beta _{pbc}$ is the same as the
coefficient of this vector in
$R^{(1)}(R^{(2)})^{-1}FS^{(1)}FB_{23}^{(1)}
B_{23}^{(2)}\beta _{pab}\otimes \beta _{pbc}$. For the first one
we have the sequence of equalities from Fig. 4.18, where the second 
equality is obtained by contracting via Theorem 4.1.
For the second one we have the equalities from Fig. 4.19, where at the
first step we used a combination of a flip and 1.b) and at the second step
we contracted. 
By moving strands around the reader can convince himself that the two are
equal.

\begin{figure}[htbp]
\centering
\leavevmode
\epsfxsize=5.1in
\epsfysize=2.8in
\epsfbox{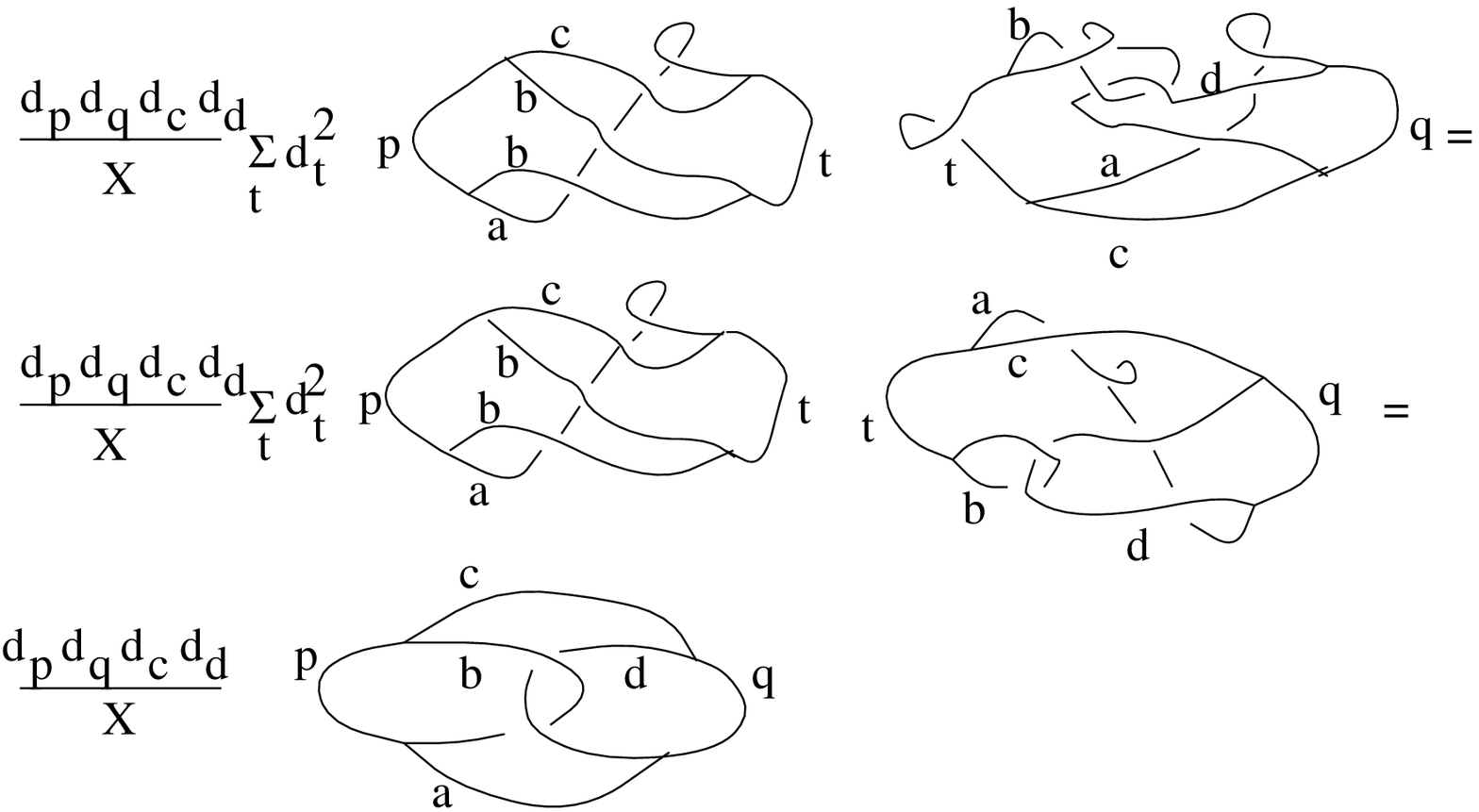}

Fig. 4.18.
\end{figure}

The groups of relations 4, 5, and 6 are straightforward. Also, we see that
the function $S$ has been chosen such that 7.a) holds. Let's prove 7.b).
Here is the place where we see why we normalized the pairing the way we
did.
We have to prove that $d_m^2d_n^2X^{-2}F\beta _{0mm}\otimes \beta _{0nn}=
\oplus _pd_md_nd_pX^{-2}\beta _{pnm}\otimes
\beta _{pmn}$. We see in Fig. 4.20 that
this is true.

\begin{figure}[htbp]
\centering
\leavevmode
\epsfxsize=5in
\epsfysize=2.8in
\epsfbox{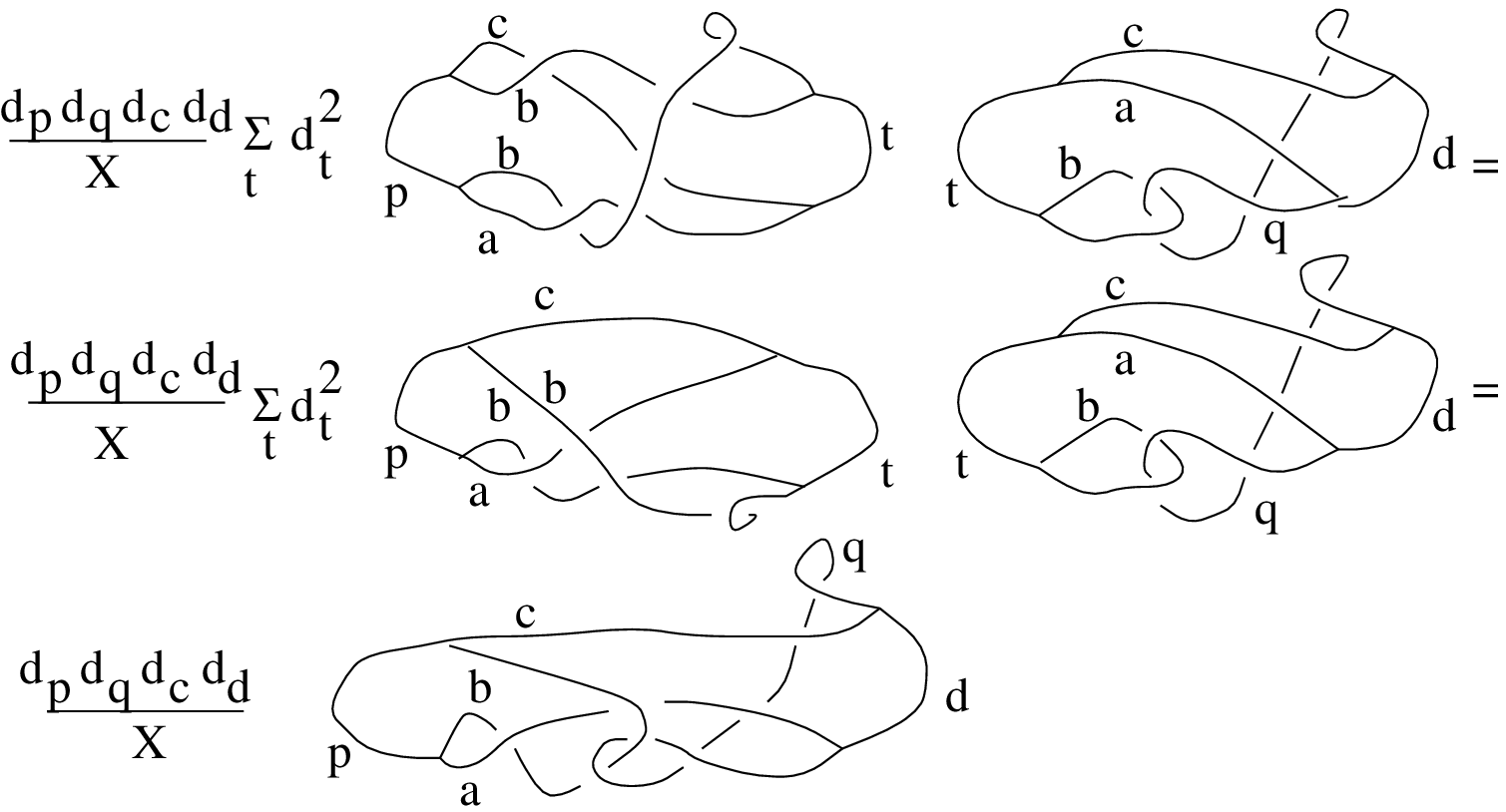}

Fig. 4.19.
\end{figure}

\vspace{15mm}

\begin{center}
{\bf 5. SOME PROPERTIES OF INVARIANTS OF 3-MANIFOLDS}
\end{center}

\bigskip

We begin this section with the generalization of Theorem 8 from [L1]
(see also Proposition 10.1 in [BHMV]) to surfaces with boundary.

\medskip

\ul{Proposition 5.1.} Let $\sa $ be a surface of genus $g$ with $n$
boundary components, and let $D$ be the DAP-decomposition of
$\sa \times S^1$ whose decomposition circles are the components of
$\partial \sa \times \{1\}$ and whose seams are of the form 
$\{x\}\times S^1$, with $x\in \partial \sa $. Then 
\begin{eqnarray*}
Z(\sa \times S^1,D,0)=\sum_{j_1,j_2,\cdots,j_n}c_{j_1,j_2,\cdots,j_n}
\beta_{j_1j_1}\otimes \beta _{j_2j_2}\otimes \cdots  \otimes \beta _{j_nj_n}
\end{eqnarray*}
where $j_1,j_2,\cdots,j_n$ run over all labelings of $\partial \sa$ and
$c_{j_1,j_2,\cdots,j_n}$ is the number of ways of labeling the diagram
in Fig. 5.1 with integers $i_k$ such that at each node we have an admissible
triple.

\begin{figure}[htbp]
\centering
\leavevmode
\epsfxsize=4.6in
\epsfysize=0.9in
\epsfbox{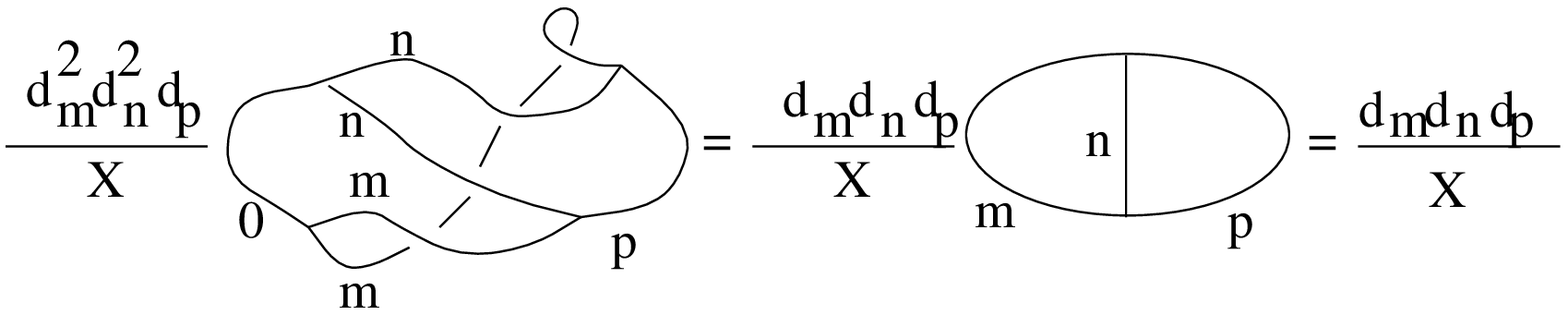}

Fig. 4.20.
\end{figure}

\begin{figure}[htbp]
\centering
\leavevmode
\epsfxsize=4.7in
\epsfysize=1.1in
\epsfbox{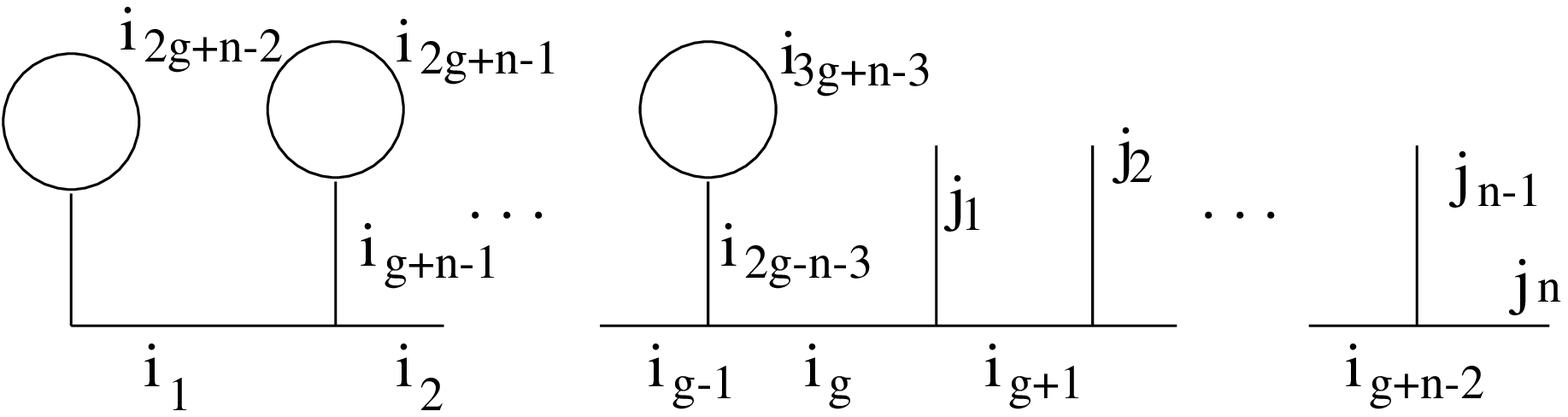}

Fig. 5.1.
\end{figure}

\ul{Proof:} Consider on $\sa $ a DAP-decomposition $D_0$ with decomposition
curves as shown in Fig. 5.2. Put on $\sa \times I$ the DAP-decomposition 
$D'$ that coincides with $D_0$ on $\sa \times \{1\}$, with $-D_0$ on
$-\sa \times \{0\}$, and on $\partial \sa \times I$ there are no extra 
decomposition circles, and the seams are vertical (i.e. of the form 
$\{ x\} \times I$).

It follows that $(\sa \times I, D',0)$ is the mapping cylinder of
$(id,0)$ (with vertical annuli no longer contracted like in the definition
of the mapping cylinder from Section 2). The mapping cylinder axiom implies 
that 
\begin{eqnarray*}
Z(\sa \times I,D',0)=\bigotimes _{j_1,j_2,\cdots,j_n}id_{j_1,j_2,\cdots,j_n}
\beta_{j_1j_1}\otimes \beta _{j_2j_2}\otimes \cdots  \otimes \beta _{j_nj_n}
\end{eqnarray*}
where $id_{j_1,j_2,\cdots,j_n}$ is the identity endomorphism on $V(\sa ,
D_0, (j_1,j_2,\cdots,j_n))$.

If we glue the ends of $\sa \times I$ via the identity map we get the 
e-3-manifold from the statement. The gluing axiom implies that in
the formula above the identity matrices get replaced by their traces.
Therefore
\begin{eqnarray*}
Z(\sa \times S^1,D,0)=\bigotimes _{j_1,j_2,\cdots,j_n}dim V(\sa ,
D_0, (j_1,j_2,\cdots,j_n))
\beta_{j_1j_1}\otimes \beta _{j_2j_2}\otimes \cdots  \otimes \beta _{j_nj_n}
\end{eqnarray*}

\begin{figure}[htbp]
\centering
\leavevmode
\epsfxsize=4.8in
\epsfysize=1.1in
\epsfbox{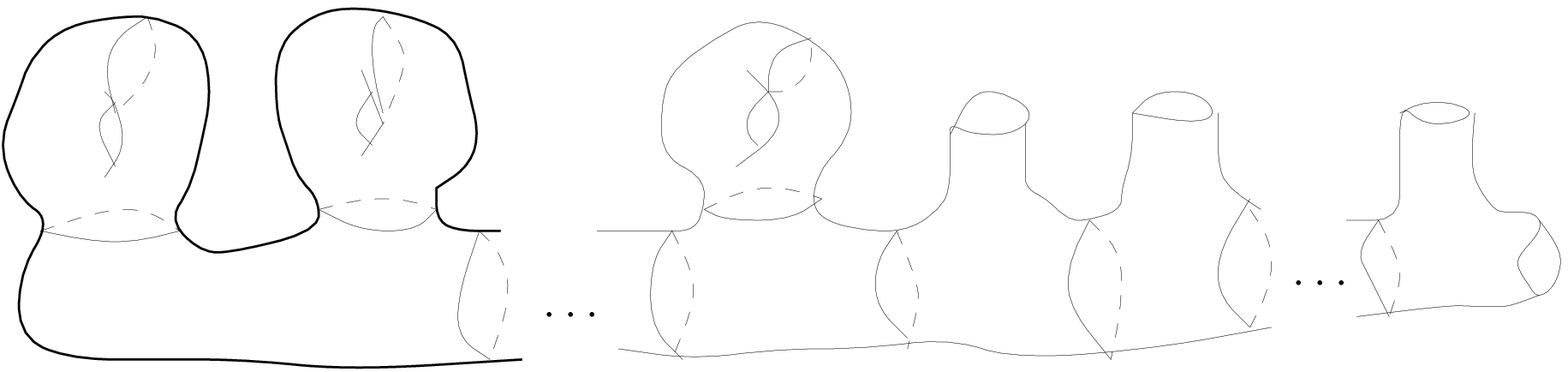}

Fig. 5.2.
\end{figure}

On the other hand the gluing axiom for $V$
implies that $dim V(\sa ,
D_0, (j_1,j_2,\cdots,j_n))
=c_{j_1,j_2,\cdots,j_n}$, which proves the 
proposition.$\Box$

\medskip

The following result shows that the Kauffman bracket not only
determines our TQFT, but also can be recovered from it. It is an
analogue of Theorem 1.1 in [G2] which showed the presence of the
skein relation of the Jones polynomial in the context of the
Reshetikhin-Turaev TQFT. Before we state the theorem we have to
introduce some notation.

Let us assume that the three e-manifolds $(M_1,D_1,0)$, $(M_2,D_2,0)$
and $(M_3,D_3,0)$ are obtained by gluing to the same e-manifold, via the same
gluing map, the
genus 2 e-handlebodies from Fig. 5.3 respectively, where the gluing
occurs along the ``exterior'' punctured spheres. Note that the three
handlebodies have the same structure on the ``exterior'' spheres, so
they produce the same change of framing (if any) when gluing.

\begin{figure}[htbp]
\centering
\leavevmode
\epsfxsize=5in
\epsfysize=1.4in
\epsfbox{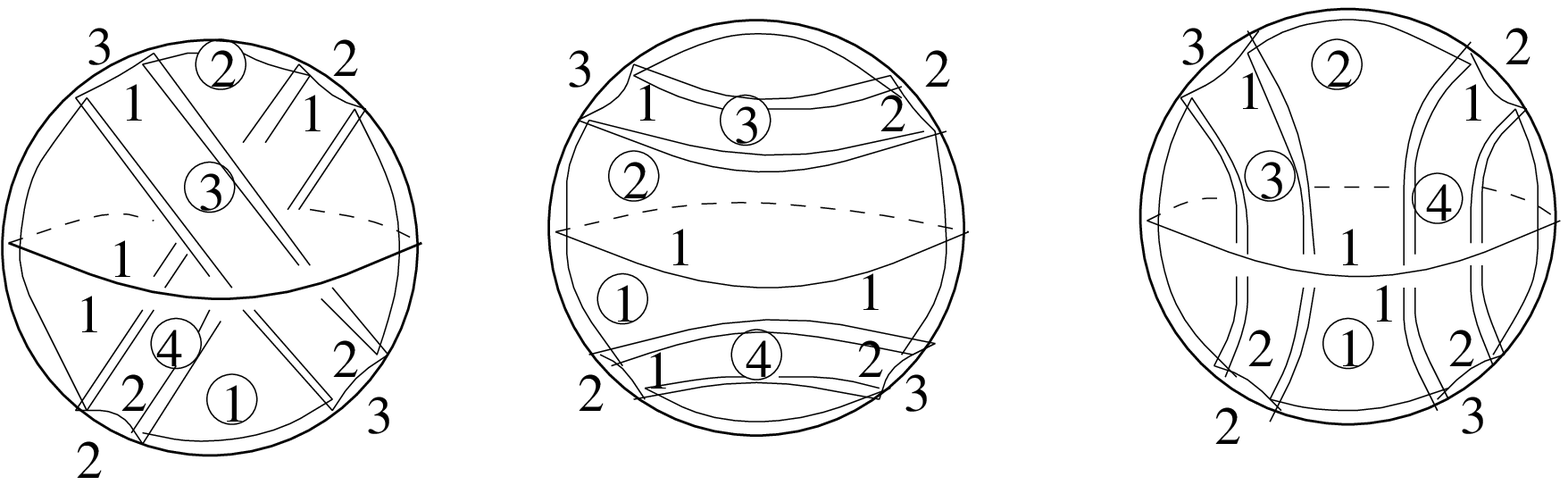}

Fig. 5.3.
\end{figure}

The ``interior'' annuli of the handlebodies are part of the boundaries of
our 3-manifolds. The gluing axiom implies that $V(\partial M_i,D_i)$ splits
as a direct sum $V_i\bigoplus V_i'$, where $V_i$ is the subspace corresponding
to the labeling of the ends of the annuli by $1$. Moreover, the gluing
axiom for $Z$ implies that $Z(M_i,D_i,0)$ also splits as $v_i\oplus
v_i'$ where $v_i\in V_i$ and $v_i'\in V_i'$. On the other hand  the spaces
$V_1$, $V_2$ and $V_3$ are canonically isomorphic. Indeed, they have a 
common part, to which the vector spaces corresponding to the two
annuli with ends labeled by $1$ are attached via the map $x\rightarrow x\otimes
\beta _{11}\otimes \beta _{11}$. Thus $v_1$, $v_2$ and $v_3$ can be thought
as lying in the same vector space. With this convention in mind,
the following result holds.

\medskip

\ul{Theorem 5.1.} The vectors $v_1$, $v_2$, and $v_3$ satisfy the
Kauffman bracket skein relation
\begin{eqnarray*}
v_1=Av_2+A^{-1}v_3.
\end{eqnarray*}

\medskip

\ul{Proof:} By the gluing axiom for $Z$ we see that it suffices to prove the
theorem in the case where $M_1$, $M_2$ and $M_3$ coincide with the 
three handlebodies (i.e. when the manifold to which they get glued is
empty).

The first e-manifold is obtained by first taking the mapping cylinder
of the homeomorphism on a pair of pants that takes the ``right leg''
over the ``left leg'' as shown in Fig. 5.4 (it should be 
distinguished from a move in the sense that it really maps one
seam into the other), then composing it with the move $B_{23}^{(1)}$,
and finally by expanding two annuli via moves of type $A^{-1}$.

\begin{figure}[htbp]
\centering
\leavevmode
\epsfxsize=2.2in
\epsfysize=0.7in
\epsfbox{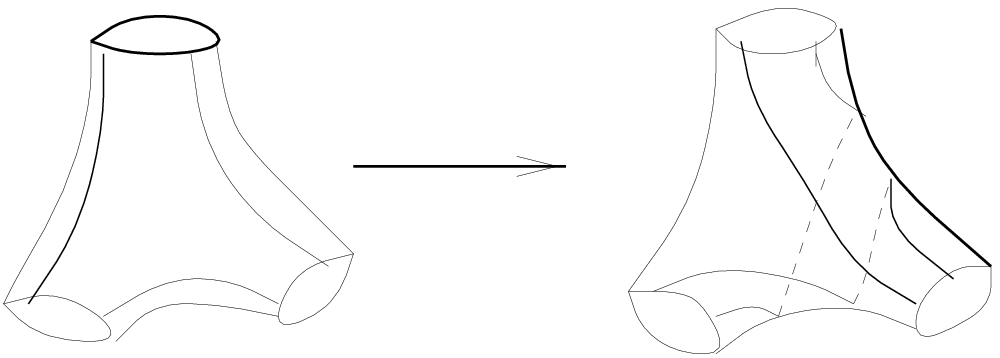}

Fig. 5.4.
\end{figure}

We get 
\begin{eqnarray*}
v_1=B_{23}\hat{\beta}_{011}\otimes \beta _{011}\otimes \beta _{11}
\otimes \beta _{11}+B_{23}\hat{\beta }_{211}\otimes \beta _{211}
\otimes \beta _{11}\otimes \beta _{11}
\end{eqnarray*}
where for $x\in V_{abc}$ we denote by $\hat{x}$ the vector in
$(V_{abc})^*$ with the property that $<x,\hat{x}>=1$.
By the definition of the pairing $\hat{\beta }_{011}=d_1^2X^{-2}\beta 
_{011}$ and $\hat{\beta }_{211}=d_1^2d_2X^{-2}\beta 
_{211}$. The computation of $B_{23}\beta _{011}$ and
$B_{23}\beta _{211}$ is described in Fig. 5.5.

\begin{figure}[htbp]
\centering
\leavevmode
\epsfxsize=6in
\epsfysize=1.3in
\epsfbox{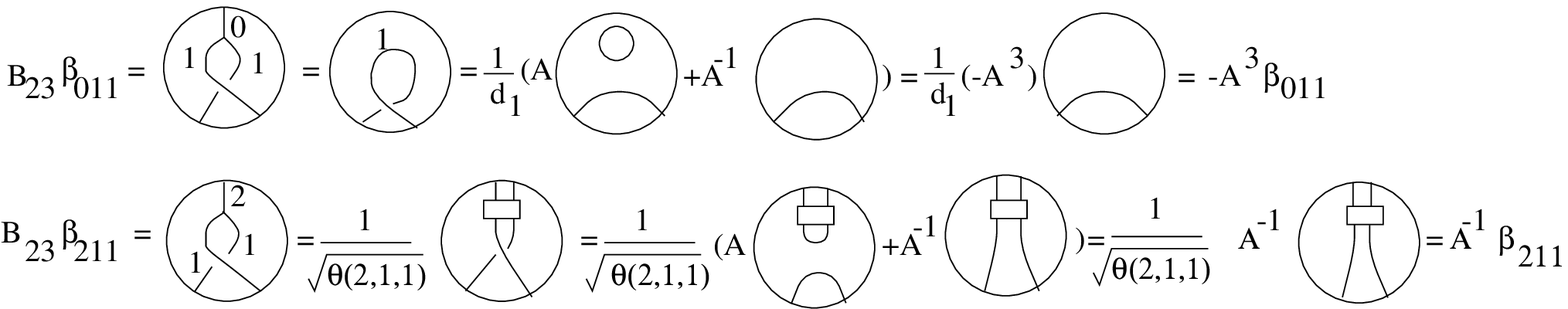}

Fig. 5.5.
\end{figure}

Hence
\begin{eqnarray*}
v_1=-A^3d_1^2X^{-2}{\beta}_{011}\otimes \beta _{011}\otimes \beta _{11}
\otimes \beta _{11}+A^{-1}d_1^2d_2X^{-2}{\beta }_{211}\otimes \beta _{211}
\otimes \beta _{11}\otimes \beta _{11}.
\end{eqnarray*}

The second manifold can be obtained by gluing along a disk the mapping
cylinders of two annuli. The mapping cylinder of an annulus has the
invariant $\oplus _a \hat{\beta }_{aa}\otimes \beta _{aa}=\oplus _a
d_a^2X^{-1}\beta _{aa} \otimes \beta _{aa}$, so after expanding a disk
and gluing the two copies together we get $\oplus _{a,b}d_a^2d_b^2X^{-2}
\beta _{0aa}\otimes \beta _{0bb} \otimes \beta _{aa}\otimes \beta _{bb}.$
But we are only interested in the component of the invariant for which 
$a=b=1$, hence $v_2=d_1^4X_{-2}\beta _{011}\otimes \beta _{011}
\otimes \beta _{11} \otimes \beta _{11}$.

Finally, the third e-manifold is the mapping cylinder of the identity
with two expanded annuli, hence
\begin{eqnarray*}
v_3=-d_1^2X^{-2}{\beta}_{011}\otimes \beta _{011}\otimes \beta _{11}
\otimes \beta _{11}+d_1^2d_2X^{-2}{\beta }_{211}\otimes \beta _{211}
\otimes \beta _{11}\otimes \beta _{11}.
\end{eqnarray*}

The conclusion follows by noting that the diagram that gives the value of
$d_1^2=\Delta _1$ is the unknot, hence $d_1^2=-A^2-A^{-2}$.$\Box$

\medskip

As a consequence of the theorem we will compute the formula for the
invariant of the complement of a regular neighborhood of a link.

\medskip

\ul{Proposition 5.2.} Let $L$ be a framed link with $k$ components,
 and $M$ be the complement of a regular neighborhood
of $L$. Consider on $\partial M$ the DAP-decomposition $D$ whose 
decomposition curves are the meridinal circles of $L$ (one for each
component) and whose seams are parallel to the framing (see Fig. 5.6.a)).
Then
\begin{eqnarray*}
Z(M,D,0)=\frac{1}{X}\sum_{n_1n_2\cdots n_k}<S_{n_1}(\alpha),
S_{n_2}(\alpha), \cdots ,S_{n_k}(\alpha )>_L \beta _{n_1n_1}
\otimes \beta _{n_2n_2}\otimes \cdots \beta _{n_kn_k}
\end{eqnarray*}
where the sum is over all labels, and $<\cdot,\cdot ,\cdots , \cdot >_L$
is the link invariant defined in Section 3.

\medskip

\ul{Proof:} We assume that $L$ is given by a diagram in the plane
with the blackboard framing.
 When $L$ is the unknot the invariant can be obtained from
Proposition 5.1 applied to the case where $\sa $ is a disk, so in this
situation $Z(M,D,0)=1/X\sa _nd_n^2\beta _{nn}$ and the formula holds.
By taking the connected sum of $k$ copies of the complement of the
unknot, and using the gluing axiom for $Z$ we see that the formula also
holds for the trivial link with $k$ components. Let us prove it in the
general case.
Put $
Z(M,D,0)={1/X}\sum_{n_1n_2\cdots n_k} c_{n_1n_2\cdots n_k} 
 \beta _{n_1n_1}
\otimes \beta _{n_2n_2}\otimes \cdots \beta _{n_kn_k}.
$ 
We want to prove that 
\begin{equation}
 c_{n_1n_2\cdots n_k}=
<S_{n_1}(\alpha),
S_{n_2}(\alpha), \cdots ,S_{n_k}(\alpha )>_L.
\end{equation}

Since by Theorem 5.1, $c_{11\cdots 1}$ and $<S_{1}(\alpha),
S_{1}(\alpha), \cdots ,S_{1}(\alpha )>_L$ satisfy both the
Kauffamn bracket skein relation, the equality holds when all
indices are equal to $1$. If some of the indices are equal to $0$, the 
corresponding link components can be neglected (by erasing them
in the case of the link, and by gluing inside solid tori in the trivial
way in the case of the 3-manifold). Therefore the equality holds if $n_i=
0,1$, $i=1,2,\cdots ,k$.

For a tuple ${\bf n}=(n_1,n_2,\cdots ,n_k)$ let $\mu ({\bf n})=
max\{n_i|i=1,2,\cdots ,k\}$ and $\nu ({\bf n})=card\{i|\ n_i=\mu ({\bf n})\}$.
We will prove (1) by induction on $(\mu ({\bf n}),\nu ({\bf n}))$, where
the pairs are ordered lexicographically. Suppose that the property is true for
all links and all tuples ${\bf n'}$ with $(\mu({\bf n'}),\nu ({\bf n'}))
< (\mu ({\bf n}),\nu ({\bf n}))$ and let us prove it for $
(\mu ({\bf n}),\nu ({\bf n}))$.

\begin{figure}[htbp]
\centering
\leavevmode
\epsfxsize=4.9in
\epsfysize=1.5in
\epsfbox{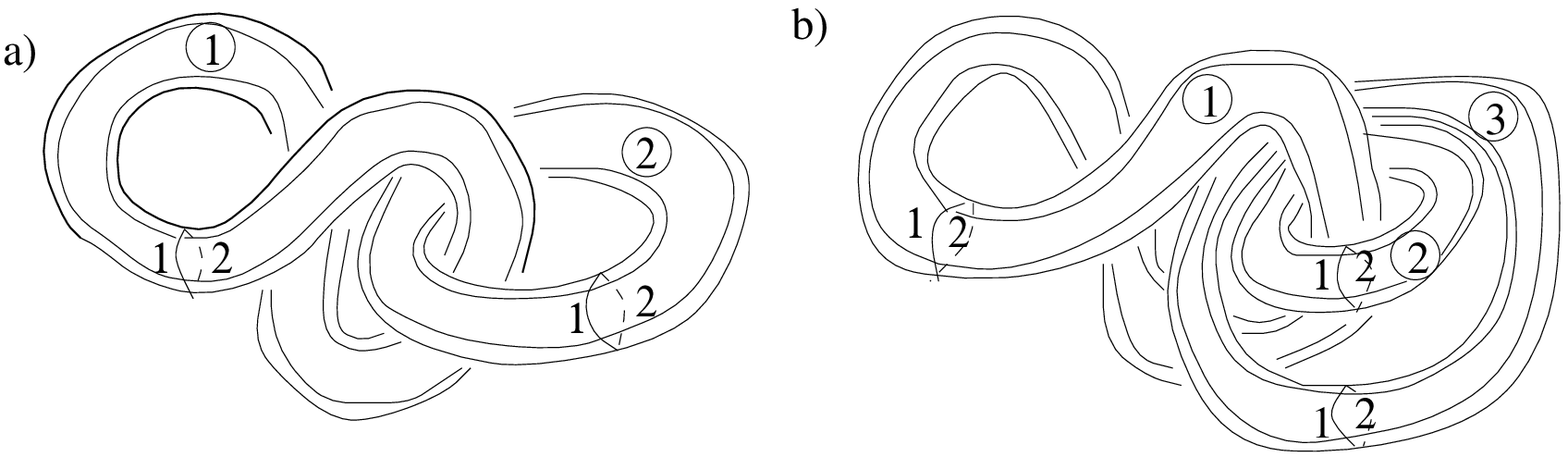}

Fig. 5.6.
\end{figure}

Let $M_0$ be the product of a pair of pants with a circle. Put on
$M_0$ a DAP-decomposition $D_0$ as described in Proposition 5.1. Then
\begin{eqnarray*}
Z(M_0,D_0,0)=\sum_{mnp}\delta _{mnp}\beta _{mm}\otimes
\beta _{nn}\otimes \beta _{pp}
\end{eqnarray*}
where $\delta _{mnp}=1$ if $(m,n,p)$ is admissible and $0$ otherwise.

Assume that in the  tuple ${\bf n}=(n_1,n_2,\cdots ,n_k)$, $n_k=
\mu ({\bf n})$. Glue $M_0$ to $M$ along the $k$-th torus of $M$ such that
in the  gluing process the DAP-decompositions of the two tori overlap.
We get an e-manifold $(M_1,D_1,0)$ that is nothing but the manifold 
associated to the link $L'$ obtained from $L$ by doubling the last
component (see Fig. 5.6. b)). 

Let $Z(M_1,D_1,0)=1/X\sa d_{m_1m_2\cdots m_k,
m_{k+1}}\beta _{m_1m_1}\otimes \beta _{m_2m_2}\otimes
\cdots \beta _{m_{k+1}m_{k+1}}$. The gluing axiom, together
with relation 6.a) from Section 3 imply that
$d_{m_1,m_2,\cdots ,m_{k+1}}=\sum_p\delta _{m_km_{k+1}p}c_{m_1,m_2,
\cdots ,m_{k-1},p}$. In particular
\begin{eqnarray*}
d_{n_1,n_2,\cdots ,n_{k-1},n_k-1,1}=
c_{n_1,n_2,\cdots ,n_k-2}+c_{n_1,n_2,\cdots ,n_k}.
\end{eqnarray*}
Applying the induction hypothesis we get
\begin{eqnarray*}
c_{n_1n_2\cdots n_k}=<S_{n_1}(\alpha ,\cdots ,S_{n_{k-1}}(\alpha ),
S_{n_k-1}(\alpha ),\alpha >_{L'}-<S_{n_1}(\alpha ,\cdots ,S_{n_{k-1}}(\alpha ),
S_{n_k-2}(\alpha )>_L.
\end{eqnarray*}

But $<S_{n_1}(\alpha ,\cdots ,S_{n_{k-1}}(\alpha ),
S_{n_k-1}(\alpha ),\alpha >_{L'}=<S_{n_1}(\alpha ,\cdots ,S_{n_{k-1}}(\alpha ),
\alpha S_{n_k-1}(\alpha )>_L$ and since $S_{n_k}(\alpha )=
\alpha S_{n_k-1}(\alpha )-S_{n_k-2}(\alpha )$ (see [L1]), we obtain the
equality in (1) and the proposition is proved.$\Box$

\medskip

\ul{Corollary} If $M$ is a closed 3-manifold obtained by performing surgery
on the framed link $L$ with $k$ components, then 
\begin{eqnarray*}
Z(M,0)=X^{-k-1}C^{-\sigma }<\omega, \omega ,\cdots , \omega >_L
\end{eqnarray*}
where $\sa $ is the signature of the linking matrix of $L$.

\medskip

\ul{Proof:} We may assume that $L$ is given by a link diagram in the
plane and its framing is the blackboard framing. Let $(M_1,D_1,0)$ be the 
e-3-manifold associated to $L$ as in the statement of 
Proposition 5.2. Consider the e-manifold $(M_2,D_2,0)$ where $M_2$ is the
solid torus and $D_2$is described in Fig. 5.7.
Applying Proposition 5.2 to the unknot we see that the invariant of this
e-manifold is $1/X\sa _nd_n^2\beta _{nn}$.

\begin{figure}[htbp]
\centering
\leavevmode
\epsfxsize=1.4in
\epsfysize=0.8in
\epsfbox{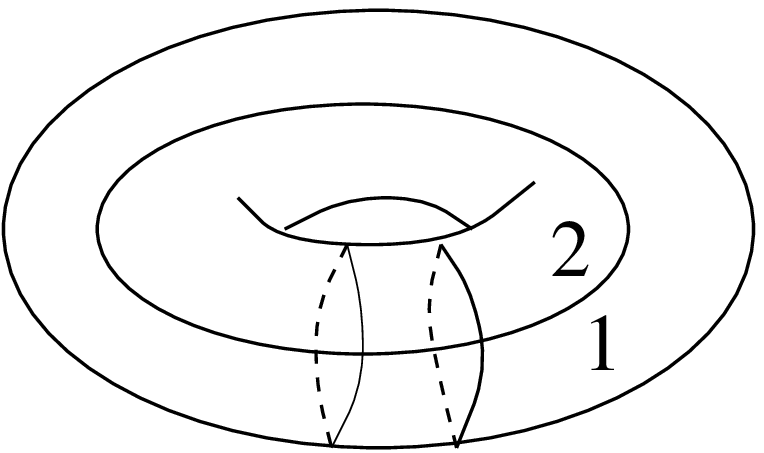}

Fig. 5.7.
\end{figure}

If we glue $k$ copies of this manifold to $M_1$ such that the
DAP-decompositions overlap we get $M$. In the gluing process the
framing changes  by $-\sigma (L_1,L_2,L_3)$ (see Section 1) where
$L_1$ is the kernel of $H_1(\partial M_1)\rightarrow H_1(M_1)$,
$L_2$ is the Lagrangian space spanned in $H_1(\partial M)$ by the
meridinal circles of the link, and $L_3$ is the one spanned by the
curves that give the framing.
It is a standard result in knot theory that $-\sigma (L_1,L_2,L_3)=
\sigma $, the linking matrix of $L$. Using the gluing axiom
for $Z$ we get
\begin{eqnarray*}
Z(M,\sigma )=X^{-k-1}\sum_{n_1,n_2,\cdots ,n_k}d_{n_1}^2d_{n_2}^2
\cdots d_{n_k}^2<S_{n_1}(\alpha ),S_{n_2}(\alpha ), \cdots ,
S_{n_k}(\alpha )>_L= X^{-k-1}<\omega, \omega ,\cdots ,\omega>_L
\end{eqnarray*}
hence
\begin{eqnarray*}
Z(M,0)=X^{-k-1}C^{-\sigma }<\omega,\omega ,\cdots ,\omega>_L.\Box
\end{eqnarray*}

\medskip

We make the final remark that this gives the invariants of 3-manifolds 
normalized as in [L1].

\medskip

{\bf REFERENCES}

\medskip

[A] Atiyah, M. F., {\em The Geometry and Physics of knots},
Lezioni Lincee, Accademia Nationale de Lincei, Cambridge Univ.
Press, 1990.

[BHMV1] Blanchet, C., Habegger, N., Masbaum, G., Vogel, P., {\em 
Topological quantum
field theories derived from the Kauffman bracket}, Topology {\bf 31}(1992),
685--699.

[Ce] Cerf, J., { \em La stratification naturelle et le th{\'{e}}oreme de la
pseudo-isotopie}, Publ. Math. I.H.E.S., {\bf 39}(1969), 5--173.

[FK] Frohman, Ch., Kania-Bartoszynska, J., {\em $SO(3)$ topological quantum
field theory}, to appear.

[G1] Gelca, R., {\em SL(2,C) topological quantum field theory with corners},
preprint, 1995.

[G2] Gelca, R., {\em The quantum invariant of the complement of a
link}, preprint, 1996.

[HT] Hatcher, A., Thurston, W., {\em A presentation of the mapping
class group of a closed orientable surface}, Topology, {\bf 19}(1980),
221--237.

[J] Jones, V., F., R., {\em Polynomial invariants of knots via von
Neumann algebras}, Bull. Amer. Math. Soc., {\bf 12}(1985), 103--111.

[K1] Kauffman, L, {\em States models and the Jones polynomial},
Topology, {\bf 26}(1987), 395--407.

[K2] Kauffman, L, {\em Knots and Physics}, World Scientific, 1991.

[Ki] Kirby, R., {\em A calculus for framed links in $S^3$}, Inventiones Math.,
{\bf 45}(1987), 35--56.

[KM] Kirby, R., Melvin, P., {\em The 3-manifold invariants of
Witten and Reshetikhin--Turaev for $sl(2,{\bf C})$}, Inventiones Math.,
{\bf 105}(1991), 547--597.

[Ko] Kohno, T., {\em Topological invariants for 3-manifolds using
representations of the mapping class groups}, Topology, {\bf 31}(1992),
203--230.

[L1] Lickorish, W.,B.,R., {\em The skein method for three-manifold
invariants}, J. Knot Theor. Ramif., {\bf 2}(1993) no. 2, 171--194.

[L2] Lickorish, W.,B.,R., {\em Skeins and handlebodies}, Pac. J. Math.,
{\bf 159}(1993), No2, 337--350.

[RT] Reshetikhin, N. Yu., Turaev, V. G., {\em Invariants of 3-manifolds
via link polynomials and quantum groups}, Inventiones Math., {\bf 103}(1991),
547--597.

[R] Roberts, J., {\em Skeins and mapping class groups}, Math. Proc. Camb.
Phil. Soc., {\bf 115}(1995), 53--77.

[T] Turaev, V., G., {\em Quantum invariants of Knots and 3-manifolds},
de Gruyter Studies in Mathematics, de Gruyter, Berlin--New York, 1994.

[Wa] Walker, K., {\em On Witten's 3-manifold invariants}, preprint, 1991.

[W] Wall, C. T. C., {\em Non-additivity of the signature}, Inventiones Math.,
{\bf 7}(1969), 269--274.
 
[We] Wenzl, H., {\em On sequences of projections}, C. R. Math.  
Rep. Acad. Sci. IX(19877), 5--9.

[Wi] Witten, E., {\em Quantum field theory and the Jones polynomial},
Comm. Math. Phys., {\bf 121}(1989), 351--399.

\medskip

Department of Mathematics, The University of Iowa, Iowa City, IA
52242 (mailing address)
{\em E-mail: rgelca@math.uiowa.edu}

Institute of Mathematics of Romanian Academy, P.O.Box 1-764,
70700 Bucharest, Romania.

\end{document}